\documentclass[reqno, 11pt]{article}
\usepackage{amsthm,amsmath,amssymb,fullpage,ifpdf}
\usepackage[usenames,dvipsnames]{color}
\usepackage[pdftex]{thumbpdf}       %create thumbnails in the PDF
\usepackage{graphics}
\usepackage{graphicx,wrapfig,floatflt,picins}   %I am using here picins, since I was not satisfied with wrapfig and floatflt. However, it is not certain which one of these 3 packages will be needed in the end. One might even resort to minipages, once the final layout is settled by the typesetter.
\usepackage{epsfig}
\usepackage{latexsym}
\usepackage{psfrag}
\graphicspath{{graphs/}}

\ifpdf
\usepackage[pdftex,letterpaper]{hyperref} %,pagebackref=true

\hypersetup{%
pdftitle={},
pdfauthor={N. Litvak, M. Vlasiou},
pdfkeywords={carousel review},
bookmarks=true,%
bookmarksnumbered=true,
pdfstartview={FitH},
pdfpagelayout={OneColumn},
colorlinks=true,
citecolor=NavyBlue,
linkcolor=RedOrange,
urlcolor=BrickRed,
bookmarksopen=true,%
bookmarksopenlevel=1,
unicode=true,%
breaklinks=true,%
}%
\else
  \newcommand\phantomsection\relax
  \newcommand{\url}[1]{#1}
  \newcommand{\href}[2]{#2}
  \usepackage{nohyperref}
\fi

\theoremstyle{plain}              %requires amsthm package
\newtheorem{theorem}{Theorem}
\newtheorem{corollary}{Corollary}
\newtheorem*{conjecture}{Conjecture}

\theoremstyle{remark}
\newtheorem{remark}{Remark}

\numberwithin{equation}{section}    %requires amsmath package

\newcommand{\e}{\mathbb{E}}
\newcommand{\p}{\mathbb{P}}

\newcommand{\suml}{\sum\limits}

\newcommand{\Dfb}[1][]{\ensuremath{F_B^{#1}}} %Memo: Df=Distribution function and df=density function
\newcommand{\Dfw}[1][]{\ensuremath{F_W^{#1}}} %\dfw(x) is normal.
\newcommand{\eDfb}[1][]{\ensuremath{\widehat{F}_B^{#1}}}
\newcommand{\eDfw}[1][]{\ensuremath{\widehat{F}_W^{#1}}}
\newcommand{\ltw}[1][]{\ensuremath{\omega^{#1}}}
\newcommand{\Dfa}[1][]{\ensuremath{F_A^{#1}}}
\newcommand{\dfw}[1][]{\ensuremath{f_W^{#1}}}
\newcommand{\lta}[1][]{\ensuremath{\alpha^{#1}}}

\begin{document}
\title{A Survey on Performance Analysis of Warehouse Carousel Systems}
\author{N.\ Litvak$^{*}$, M.\ Vlasiou$^{**}$}
\date{\today}
%\keywords{}
\maketitle

\begin{center}
$^{*}$ Faculty of Electrical Engineering, Mathematics and Computer Science,
\\      Department of Applied Mathematics, University of Twente,
\\      7500 AE Enschede, The Netherlands.
\end{center}

\begin{center}
$^{**}$ Eurandom and Department of Mathematics \& Computer Science,
\\ Eindhoven University of Technology,
\\ P.O.\ Box 513, 5600 MB Eindhoven, The Netherlands.
\\ \vspace{0.3cm}\href{mailto:n.litvak@ewi.utwente.nl}{n.litvak@ewi.utwente.nl}, \href{mailto:m.vlasiou@tue.nl}{m.vlasiou@tue.nl}
\end{center}

\begin{abstract}
This paper gives an overview of recent research on the performance evaluation and design of carousel systems. We  discuss picking strategies for problems involving one carousel, consider the throughput of the system for problems involving two carousels, give an overview of related problems in this area, and present an extensive literature review. Emphasis has been given on future research directions in this area.

%The focus is on results, intuition, and insight rather than
%methods and techniques.
\end{abstract}
\medskip

\medskip


\section{Introduction}\label{s:intro}
A carousel is an automated storage and retrieval system, widely
used in modern warehouses. It consists of a number of shelves or
drawers, which are linked together and are rotating in a closed
loop. It is operated by a picker (human or robotic) that has a
fixed position in front of the carousel. A typical vertical
carousel is given in Figure~\ref{fig:carousel}.

Carousels are widely used for storage and retrieval of small and
medium-sized items, such as health and beauty products, repair
parts of boilers for space heating, parts of vacuum cleaners and
sewing machines, books, shoes and many other goods. In e-commerce
companies use carousel to store small items and manage small
individual orders. An {order} is defined as a set of items that
must be picked together (for instance, for a single customer).

%\begin{floatingfigure}[p]{7 cm}
%\scalebox{0.3}{\includegraphics*[1,5][610,790]{carousel}}
%\caption{A typical vertical carousel}
%\noindent \hrulefill
% %\label{fig:carousel}
%\end{floatingfigure}

%
%\begin{figure}
%\begin{center}
%\includegraphics[width=0.5\textwidth]{carousel}
%\end{center}
%\caption{A typical vertical carousel}
%  %\label{fig:carousel}
%\end{figure}

Carousels are highly versatile, and come in a huge variety of configurations, sizes, and
types. They can be horizontal or vertical and rotate in either one
or both directions. Although both {unidirectional} (one-way
rotating) or {bidirectional} (two-way rotating) carousels are
encountered in practice, the bidirectional types are the most
common (as well as being the most efficient)~\cite{hassini03}. One
of the main advantages of carousels is that, rather than having
the picker travel to an item (as is the case in a warehouse where
items are stored on shelves), the carousel rotates the items to
the picker. While the carousel is travelling, the picker has the
time to perform other tasks, such as pack or label the retrieved
items, or serve another carousel. This practice enhances the
operational efficiency of the warehouse.

Carousel models have received much attention in the literature and
continue to pose interesting problems. There is a rich literature
on carousels that dates back to 1980~\cite{weiss80}. In
Section~\ref{s:literature} we shall review some of the main
research topics that have been of interest to the research
community so far. To name a few, one may wish to study various
ways of storing the items on a carousel (storage arrangements) so as to minimise the total time needed until an order is completed (response
time) or the strategy that should be followed in rotating the
carousel so as the total time the carousel travels between items
of one order is minimised (travel time for a single order). One
may also consider design issues, for instance, the problem of
pre-positioning the carousel in anticipation of storage or
retrieval requests (choosing a \textit{dwell point}) in order to
improve the average response time of the system. The list of
references presented here is by no means exhaustive; it rather
serves the purpose of indicating the continuing interest in
carousels.
\piccaption{A typical vertical carousel.\label{fig:carousel}}
\parpic[l]{\includegraphics[width=0.4\textwidth]{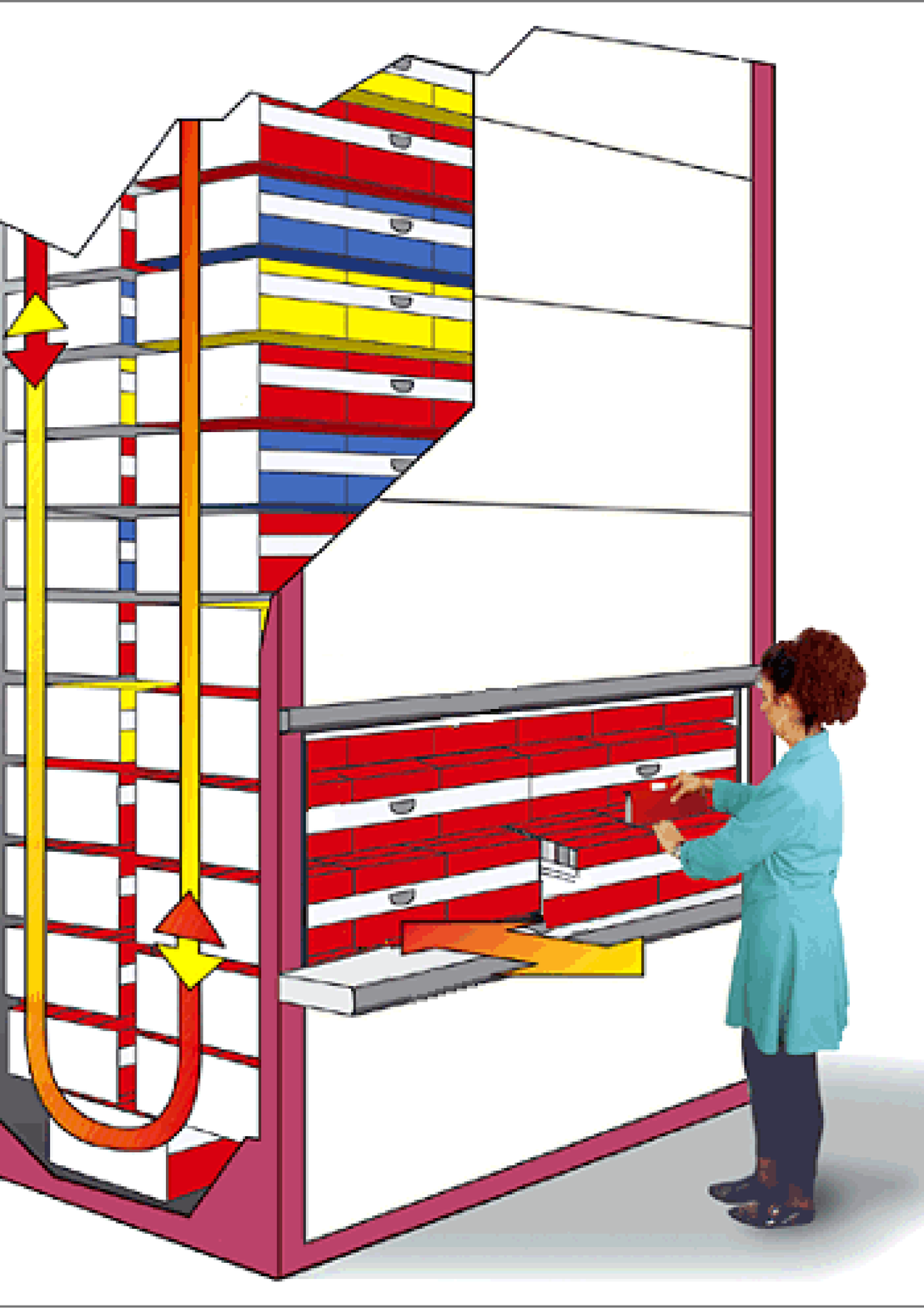}}

%note to typesetter: I had issues with the packages wrapfig, floatfig, floatflt. the package picins used with parpic will put the picture at the beginning of a paragraph. it comes with all sort of options for the format and positioning which may help in the final formatting.

In this review paper we focus on the modelling and the performance
of carousel systems. Usually a carousel is modelled as a circle,
either as a discrete model
\cite{bartholdi86,jacobs00,stern86,yeh02}, where the circle
consists of a fixed number of locations, or as a continuous one
\cite{ghosh92,litvak,vdberg96,vlasiou04}, where the circle has
unit length and the locations of the required items are
represented as arbitrary points on the circle. Throughout this
paper we shall view the carousel as a continuous loop of unit
length. Beyond this initial assumption, we shall examine modelling
issues such as how to model travel times or picking times of items
in a system of several carousels so as to be able to derive
approximations of various performance characteristics. Under
``performance'' one may understand a variety of notions. For
example, in single-carousel single-order problems (cf.\
Section~\ref{s:nelly}), the performance measure under
consideration is the travel time of the carousel until all items
in an order are picked. On the other hand, in
Section~\ref{s:maria}, performance may be measured by the time the
picker is idle between picking items from various carousels, i.e.\
by the picker's utilisation.

In this paper we consider two research topics in detail. In
Section~\ref{s:nelly}, we discuss the problem of choosing a
reasonable picking strategy for one order and a single carousel,
where the order is represented as a list of items, and by order
pick strategy we mean an algorithm that prescribes in which
sequence the items are to be retrieved. We present a general
probabilistic approach developed by Litvak \emph{et
al}.~\cite{litvak,litvak01,litvak02,litvak01a} to analytically
derive the probability distribution of the travel time in case
when items locations are independent and uniformly distributed.
This line of research seems to be the only example in the
literature where exact statistical characteristics of the travel
time have been obtained by means of a systematic mathematical
approach. The presented technique is based on properties of
uniform spacings and their relations to exponential distributions.
We demonstrate the effectiveness of this method by considering
several relevant order-picking strategies, such as the greedy
nearest-item strategies and so-called $m$-step strategies that
provide a good approximation for the optimal (shortest) route.

In Section~\ref{s:maria} we consider the second topic that relates
to multiple-carousel settings and the modelling challenges that
appear in such problems. Having optimised the travel time of a
single carousel for a single order, one wonders if optimising
locally every time each order on each carousel leads to the best
solution (fastest, cheapest, or with the largest picker
utilisation) for a complicated system. As is mentioned later on,
multiple-carousel problems become too complicated too quickly, and
often exact analysis is not possible. Therefore, we discuss which
concessions have to be made in order to be able to obtain
estimates of the performance measures we are interested in, and we
give in detail the impact that these concessions have on our
estimations. There exist a few exact results for two-carousel
models and related models in healthcare logistics; see Boxma and
Vlasiou \cite{boxma07} and Vlasiou {\em et al.}~\cite{vlasiou07a}--\cite{vlasiou07}. However, to the best of
our knowledge, no exact results exist for systems involving more
than two carousels.

Preferably, these two research topics that we consider in this
paper should be studied in parallel. However, establishing any
exact results, say on determining the optimal retrieval and
travelling strategy for a multiple-carousel model, without any
restrictions to the sequence the items in an order are picked or
the sequence the carousels are served, seems to be intractable.
Nonetheless, quite a few research opportunities related to the
optimal design and control of carousel systems are still
available. We elaborate on further research topics in
Section~\ref{s:research}. We conclude with
Section~\ref{s:literature}, which outlines the problems examined
so far on carousels and related storage and retrieval systems.

\section{Picking a single order on a single carousel}\label{s:nelly}

Performance analysis of single units is a necessary step in
structural design of order pick systems~\cite{yoon96}. In a
setting of  a single order on a single carousel, the major
performance characteristic is the response time, that is, the
total time it takes to retrieve an order. The response time
consists of pick times needed to collect the items from their
locations by an operator, and the travel (rotation) time of the
carousel.  While pick times can hardly be improved, the travel
time depends on the location of each item and the order picking
sequence, and thus, it is subject to analysis and optimisation.
Therefore, in this section, we discuss properties of the travel
time needed to collect an order of $n$ items. In this section, our
focus is on the case when the item locations are randomly
distributed on a carousel circumference. This model allows one to
compute statistical characteristics of the travel time such as the
average travel time or the travel time distribution. Later on, in
Section~\ref{ss:Picking a single order} we discuss some results
from the literature on evaluating the travel times under different
assumptions on the items locations, in particular, the case when
the pick positions are fixed.

We note that in case of a single carousel, it is natural to assume
that the pick times and the travel time are independent. The
situation, however, is quite different in the systems of two or
more carousels, where pick times on one carousel affect the travel
times on other carousels. This issue will be discussed in detail
in Section~\ref{s:maria}.

The model addressed in this section is as follows. We model a
carousel as a circle of length~1. The order is represented by the
list of $n$ items whose positions are independent and uniformly
distributed on $[0,1)$. For ease of presentation, we act as if the
picker travels to the pick positions instead of the other way
around. Also, we assume that the acceleration/deceleration time of
the carousel is negligible or that it is assigned to the pick
time, and that the carousel rotates at unit speed. Therefore
the travel distance can be identified with the travel time (see
also Section~\ref{ss:Design issues}).

\begin{figure}[htb]
%\begin{figure}[e]
\begin{center}
\includegraphics[width=1.5in]{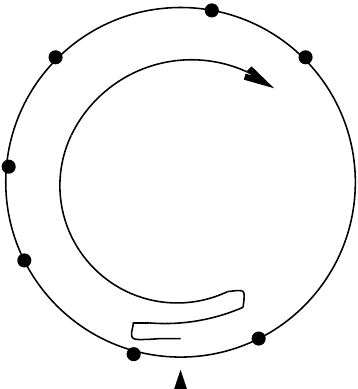}
\end{center}
%\centering {\epsfxsize=1.5in \epsfbox{ni_route.eps}}
\caption{{\small A route under the nearest-item heuristic. \label{fig_ni_route}}}
\end{figure}

Obviously, the travel time depends heavily on the pick strategy.
Here by {\it order pick strategy} we mean an algorithm that
prescribes the sequence in which the items are collected. For
example, assume that the picker always proceeds in the clockwise
(CW) direction and denote by $T_n^{CW}$ the time needed to collect
$n$ items under this simple strategy. Then, clearly, the
distribution function $\mathbb{P}(T_n^{CW}\le t)$ of $T_n^{CW}$
simply equals $t^n$, $0<t\le 1$. However, we would like to study
strategies that provide smaller travel times. In this sense, a
better algorithm that one can think of is the `greedy' strategy,
also called the {\it nearest-item heuristic}: always travel to the
nearest item to be picked (as in Figure~\ref{fig_ni_route}).
The nearest-item strategy indeed performs very well and is often
used in practice, but the question is: ``what is the distribution
of the travel time under the nearest-item heuristic?''. This
problem is not at all trivial. For example, straightforward
methods, such as conditioning on possible item locations, do not
lead to feasible calculations. The same applies to the optimal
strategy. Bartholdi and Platzman~\cite{bartholdi86} showed that
the shortest route admits at most one turn. Intuitively, this
follows merely by observing Figure~\ref{fig_ni_route}, where the
displayed route can be shortened by collecting the first item in
the counterclockwise direction and then collecting the rest of the
items rotating clockwise. Thus, the shortest route is merely the
minimum among the $2n$ candidate routes than have at most one
turn. However, in spite of this simple structure of the shortest
route, its distribution function is hard to derive.

Below we discuss in detail a general methodology developed by
Litvak {\em et al.}~\cite{litvak,litvak01,litvak02,litvak01a} to
obtain the distribution of the travel time under various order
pick strategies. The proposed technique is based on properties of
uniform spacings and their connection with exponential random
variables. We show how this approach allows us to derive exact and
often counterintuitive results on several relevant order pick
strategies. Some other methods from the literature are described
in Section~\ref{ss:Picking a single order}.

We start with introducing the notation and presenting some background
results. Let the random variable $U_0=0$ be the picker's starting
point and the random variable $U_i$, where $i = 1, 2, \ldots, n$,
be the position of the $i$th item. We suppose that the $U_i$'s, $i
= 1, 2, \ldots, n$, are independent and uniformly distributed on
$[0,1)$. Let $U_{1:n}, U_{2:n},\ldots U_{n:n}$ denote the order
statistics of $U_{1}, U_{2},\ldots U_{n}$ and set $U_{0:n}=0$,
$U_{n+1:n}=1$. Then the uniform spacings are defined as
\begin{equation}
\label{eq_spacings} D_{i,n}=U_{i:n}-U_{i-1:n},\quad 1\le i\le n+1.
\end{equation}
If we consider $n$ items randomly located on a circle, then the
spacings $D_{2,n}, D_{3,n}, \ldots, D_{n,n}$ are the distances
between two neighbouring items, and the spacings $D_{1,n}$ and
$D_{n+1,n}$ are the distances between the starting point and the
two items adjacent to it. Whatever strategy the picker takes, he
always has to cover one or more uniform spacings on his way from
one location to another. Hence, in general, the travel time can be
expressed as a function of the uniform spacings.

Uniform spacings have been analysed extensively in two classical
review papers by Pyke~\cite{pyke65,pyke72}. The author gives four useful
constructions that establish a connection between uniform spacings
and exponential random variables. We will use such a connection in
the following form. Let $X_1$, $X_2,\ldots$ be independent
exponential random variables with mean~1. Moreover, define the random variables
\[S_0=0, \qquad S_i=X_1+X_2+\cdots+X_i,\quad i\ge 1.\] Then,
according to Pyke~\cite{pyke65}, uniform spacings can be
represented as follows:
\begin{equation}
\label{eq_prop} (D_{1,n},D_{2,n},\ldots,
D_{n+1,n})\stackrel{d}{=}\left(X_1/S_{n+1}, X_2/S_{n+1}, \ldots,
X_{n+1}/S_{n+1}\right).
\end{equation}
Here and throughout this paper $a\stackrel{d}{=} b$ means that $a$
and $b$ have the same probability distribution. Linear
combinations of uniform spacing have nice properties. In
particular, the moments of linear combinations with non-negative
coefficients can be easily computed, and their distribution
function has been derived by Ali~\cite{ali73}, Ali and
Obaidullah~\cite{ali82}.

Now, let $X$ and $Y$ be independent exponential random variables
with parameters $\lambda$ and $\mu$, respectively. We write
$X=X_1/\lambda$, $Y=Y_1/\mu$, where $X_1$ and $Y_1$ are
independent exponential random variables with parameter~1. Then,
given the event $[X<Y]$, we obtain the following useful
statements: (i)~the distribution of $X=\min\{X,Y\}$ is exponential
with parameter $\lambda+\mu$ (property of the minimum of two
exponentials), which is distributed as $X_1/(\lambda+\mu)$;
(ii)~since $[Y>X]$, then, according to the memoryless property,
$Y$ can be written as a sum of two terms: $\min\{X,Y\}$ and
another independent exponential with parameter $\mu$, so $Y$ is
distributed as $X_1/(\lambda+\mu)+Y_1/\mu$. (iii)~it is easy to
check that the distribution of $S=\lambda X+\mu Y=X_1+Y_1$ is
independent of the event $[X<Y]$ because according to (i) and
(ii), given $[X<Y]$, $S$ is again distributed as $X_1+Y_1$ (see
also Chapter~2 of \cite{litvak}).

Based on the above-mentioned properties of exponential random
variables, and their connections to uniform spacings and travel
times, one may adopt the following methodology for analysing the
travel times under various
strategies~\cite{litvak,litvak01,litvak02,litvak01a}:
\begin{enumerate}
\item Express the travel time under a given strategy as a function
of uniform spacings. \item By conditioning on linear inequalities
between the spacings and employing the above mentioned properties
of exponential random variables, rewrite the travel time as a
linear combination of uniform spacings or as a probabilistic
mixture of such linear combinations. \item Use the results
from~\cite{ali73,ali82} to obtain the moments and the
distribution of the travel time.
\end{enumerate}
Below we show how this approach works in case of the nearest-item
heuristic~\cite{litvak01,litvak01a} and so-called $m$-step
strategies~\cite{litvak02}.

\subsection{The nearest-item heuristic}\label{ss:narest-item}
\begin{figure}[htb]
\begin{center}
\includegraphics[width=9cm]{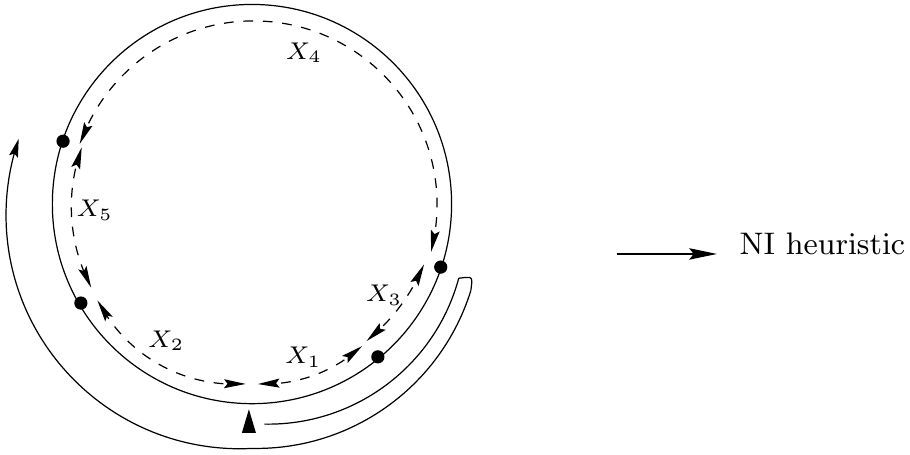}
\caption{The nearest-item route of the picker facing 5 exponential spacings.\label{fig_travel}}
\end{center}
\end{figure}

Under the nearest-item heuristic, the picker always moves towards
the nearest item to be retrieved. The positions of the items
partition the circle in $n+1$ uniform spacings $D_{1,n},D_{2,n}
,\ldots,D_{n+1,n}$ defined by \eqref{eq_spacings}. Under the
nearest-item heuristic, the picker  first considers the two
spacings adjacent to his starting position and then travels to the
nearest item. Next he also looks at the other spacing adjacent to
that item and compares the distance to the item located at the
endpoint of that spacing and the distance to the first item in the
other direction, which is the sum of the spacings previously
considered. Then he travels again to the nearest item, and so on.
Furthermore, by employing~\eqref{eq_prop}, we may act as if the
picker faces non-normalised exponential spacings, and afterwards
divide the travel time (which is equal to the travel distance) by
the sum of all spacings. Then it is clear that each new spacing
faced by the picker is independent of the ones already observed.
Now let $X_i$, $i = 1, \ldots, n+1$, denote the $i$-th
non-normalised exponential spacing faced by the picker. That is,
the spacings are numbered as observed by the picker operating
under the nearest-item heuristic (see Figure~\ref{fig_travel}).
%
%\input{ni9_fig.tex}
%
%\begin{figure}[htb]
%\centering
%%
%\psfrag{NI heuristic}[][]{{\small\hspace{.2cm} NI heuristic}}
%\psfrag{X1}[][]{$\scriptstyle{X_1}$}
%\psfrag{X2}[][]{$\scriptstyle{X_2}$}
%\psfrag{X3}[][]{$\scriptstyle{X_3}$}
%\psfrag{X4}[][]{$\scriptstyle{X_4}$}
%\psfrag{X5}[][]{$\scriptstyle{X_5}$}
%\includegraphics[width=9cm]{ni9}
%\caption{\small The nearest-item route of the picker facing 5
%exponential spacings.\label{fig_travel}}
%\end{figure}
%
Then $T^{NI}_n$ can be expressed as
\begin{equation}
\label{eq_ni_tt} T^{NI}_n = \sum_{i=2}^{n+1} \frac{\min\{X_i,
S_{i-1}\}}{S_{n+1}}.
\end{equation}

We first provide an informal explanation of how the proposed
methodology can be applied to \eqref{eq_ni_tt}. To start with, note that
first term in the right-hand side of \eqref{eq_ni_tt} is
$\min\{X_1,X_2\}/S_{n+1}$, which is distributed simply as
$(1/2)X_1/S_{n+1}$. Moreover, under the event $[X_1<X_2]$ the rest
of the sum remains unaltered. Further, consider the term
\begin{equation}
\label{eq_2ndterm}(1/2)X_1+\min\{X_3,S_2\}=(1/2)X_1+\min\{X_3,X_1+X_2\}.
\end{equation}
Let $X'_1$, $X'_2$, $X'_3$ be auxiliary independent exponential
random variables with mean~1. Given $[X_3<X_1]$, the random
variable $X_3$ is distributed as $(1/2)X'_1$, $X_1$ is distributed
as $(1/2)X'_1+X'_2$ and $X_2$ is distributed as $X'_3$. Then the
term in \eqref{eq_2ndterm} is distributed as
$(3/4)X'_1+(1/2)X'_2$. Furthermore, given the event $[X_3>X_1,
X_3<X_1+X_2]$, we obtain that $X_1$ is distributed as $(1/2)X'_1$,
$X_3$ is distributed as $(1/2)X'_1+(1/2)X'_2$ and $X_2$ is
distributed as $(1/2)X'_2+X'_3$. Substituting the above in
\eqref{eq_2ndterm}, we obtain again $(3/4)X'_1+(1/2)X'_2$!
Remarkably, under the event $[X_3>X_1+X_2]$, \eqref{eq_2ndterm}
again transforms into $(3/4)X'_1+(1/2)X'_2$. Furthermore, the sum
$S_3=X_1+X_2+X_3$ becomes simply $S_3=X'_1+X'_2+X'_3$. We may now
rename $(X'_1,X'_2,X'_3)$ back to $(X_1,X_2,X_3)$ since the two
3-dimensional vectors are identically distributed. Then the term
\eqref{eq_2ndterm} becomes $(3/4)X_1+(1/2)X_2$, and the rest of
the terms in the right-hand side of \eqref{eq_ni_tt} remain
unaltered in all three cases. Proceeding further, we obtain the
next statement which is proved rigorously in~\cite{litvak01}.
\begin{theorem}
[Litvak and Adan~\cite{litvak01}] \label{th_nearest_item} For all
$n=1,2,\ldots$,
\begin{equation}
\label{eq_ni_repr}T^{NI}_n\stackrel{d}{=}
\sum_{i=1}^n\left(1-\frac{1}{2^i}\right)D_{i,n}
\end{equation}
and
\begin{equation}\label{eq_ni_dist}
{\mathbb{P}}(T^{NI}_n\le t) = \sum_{i=0}^{n} \left(2^it-2^i+1\right)_+^n\prod_{j=0\atop j\ne i
}^n\frac{2^j}{2^j-2^i},\quad 0<t\le 1,
\end{equation}
where $x_+=x$ if $x>0$ and $x_+=0$ otherwise.
\end{theorem}
Here \eqref{eq_ni_dist} follows directly from \eqref{eq_ni_repr}
and the result by Ali~\cite{ali73}, which we applied in the form
given by Theorem~2 in \cite{ali82}.

The above theorem is surprising because it provides an elegant
solution for a problem that looks intractable at first. An interesting by-product is the distribution of the number
of turns under the nearest-item heuristics and the
counterintuitive result that {\it the travel time and the number
of turns are independent}~\cite{litvak}! The latter can be
seen directly from \eqref{eq_ni_tt}. Indeed, a turn after step $i$
is equivalent to the event $[X_{i+1}>S_i]$. However, as we saw
earlier, the form of the distribution of the travel time is given by \eqref{eq_ni_repr} and it is independent of this sort of events.

\subsection{The m-step strategy}\label{ss:m-step}
\begin{figure}[htb]
\begin{center}
\includegraphics[width=9cm]{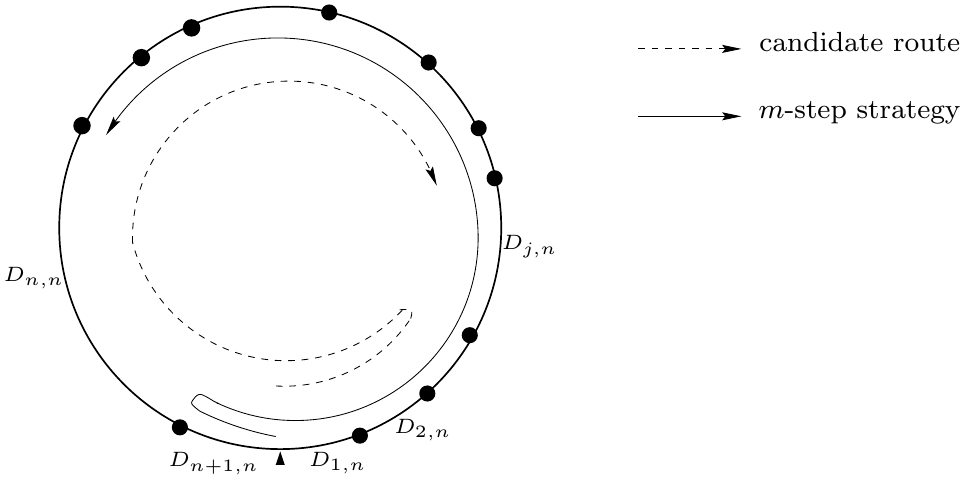}
\caption{A route under the $m$-step strategy.\label{fig:m-step}}
\end{center}
\end{figure}

Under the $m$-step strategies, the picker chooses the shortest
route among the $2(m+1)$ routes that change direction at most
once, and only do so after collecting no more than $m$ items. Note
that the optimal strategy is in fact an $(n-1)$-step strategy
since it is never optimal to turn more than once, and the maximal
possible number of items collected before a turn is $n-1$. The
$m$-step strategies give a good approximation for the shortest
travel time. In fact,  they often provide the optimal route even
for moderate values of $m$, as in Figure~\ref{fig:m-step}.
Rouwenhorst {\em et al.}~\cite{rouwenhorst96} were the first to
propose these strategies as an upper bound for the optimal route.
In case of independent uniformly distributed pick positions, they
obtained the distribution of the travel time under the $m$-step
strategy for $m\le 2$ using analytical methods. Later on, Litvak
and Adan~\cite{litvak02} applied the described methodology based
on the properties of uniform spacings to completely analyse the
travel time under the $m$-step strategies, provided $2m<n$. The
travel under the $m$-step strategy can be expressed as follows
$$
T^{(m)}_{n}=1-\max\left\{\max_{1\le j\le
m+1}\left\{D_{j,n}-\sum_{l=1}^{j-1}D_{l,n}\right\}, \max_{1\le
j\le
m+1}\left\{D_{n+2-j,n}-\sum_{l=1}^{j-1}D_{n+2-l,n}\right\}\right\}.
$$
Indeed, the term $D_{j,n}-\sum_{l=1}^{j-1}D_{l,n}$ is the gain in
travel time (compared to one full rotation) obtained by skipping
the spacing $D_{j,n}$ and going back instead, ending in a
clockwise direction. On the other hand,
$D_{n+2-j,n}-\sum_{l=1}^{j-1}D_{n+2-l,n}$ is the gain obtained by
skipping the spacing $D_{n+2-j,n}$ and going back ending
counterclockwise. Under the $m$-step strategy the picker skips the
spacing that provides the largest possible gain (see
Figure~\ref{fig:m-step}).
%
%\input{optimal_fig}
%
%\begin{figure}[htb]
%%\label{fig1}
%\centering {\footnotesize
%%
%\psfrag{m}[][]{$$}
%\psfrag{n+1-m}[][]{$$}
%\psfrag{D1}[][]{$\scriptstyle{D_{1,n}}$}
%\psfrag{D2}[][]{$\scriptstyle{D_{2,n}}$}
%\psfrag{Dm+1}[][]{
%$$}
%\psfrag{Dn+1-m}[][]{
%$$} \psfrag{Dj}[][]{$\scriptstyle{D_{j,n}}$}
%\psfrag{D_{m+1}}[][]{$$}
%\psfrag{Dn}[][]{\hspace*{-.2cm}$\scriptstyle{D_{n,n}}$}\psfrag{Dn+1}[][]
%{$\scriptstyle{D_{n+1,n}}$}
%\psfrag{ldots}[][]{$\scriptstyle{\ldots}$} \psfrag{candidate route}[][]{\hspace{.4cm} candidate route}
%%
%\psfrag{m-step strategy}[][]{\hspace*{.4cm} $m$-step strategy}
%%
%\includegraphics[width=9cm]{fig1_orl}
%%
%\caption{A route under the $m$-step strategy.\label{fig:m-step}}}
%%
%\end{figure}
%
Using property \eqref{eq_prop}, and after appropriate
manipulations of exponential random variables, one can prove the
following result.

\begin{theorem}[Litvak and Adan~\cite{litvak02}]
\label{th_m-step} For any $m=0,1,\ldots$, with $2m<n$,
\begin{equation}\label{eq_m-step_repr}
T^{(m)}_{n}\stackrel{d}{=}1-\frac{1}{S_{n+1}}
\max\left\{\sum_{j=1}^{m+1}\frac{1}{2^{j}-1}\,X_{j},
\sum_{j=1}^{m+1}\frac{1}{2^{j}-1}\,{X_{n+2-j}}\right\}.
\end{equation}
\end{theorem}

The maximum in the right-hand side of \eqref{eq_m-step_repr}
implies that $T_n^{(m)}$ is distributed as a complicated
probabilistic mixture of linear combinations of uniform
spacings~\cite{litvak02}. The number of terms in this mixture is
the well-known Catalan number
\[\frac{1}{m+2}{{2m+2}\choose{m+1}},\]
which grows extremely fast with $m$. Computing the expectations,
we conclude that on average, the $m$-step strategy performs better
than the nearest-item heuristic already for $m=2$ provided $n\ge
5$.

Again, as a by-product, we can obtain the distribution of the
number of steps before the turn. Moreover, the latter random
variable turns out to be independent of the travel time. This
surprising statement follows from a similar reasoning as the
independence of the travel time and the number of turns under the
nearest-item heuristic. Furthermore, when $n$ goes to infinity,
the number of steps before the turn converges to a shifted
geometric distribution with parameter $1/2$. That is, in the
limit, with probability $1/2$ there will be no turn, with
probability $1/4$ there will be one step before a turn, etc. Also,
in the limit, the $m$-step strategy with $2m<n$ coincides with the
optimal strategy since the probability of achieving the minimal
travel time by making more than $n/2$ steps before a turn will
converge to zero. Thus, for large enough $n$, the probability that
a 2-step strategy provides an optimal route is about $7/8$. This
explains the remarkably good performance of the $m$-step
strategies.

As a side remark, we would like to note that \cite{litvak01b} provides slightly more general results than those presented in
\eqref{eq_ni_repr} and \eqref{eq_m-step_repr}.

\subsection{Optimal route}\label{ss:optimal strategy}

Since the optimal strategy simply coincides with the $(n-1)$-step
strategy (at most one turn after collecting at most $n-1$ items)
it can be analysed by methods from
Section~\ref{ss:m-step}. However, the condition $2m<n$ is violated
for $m=n-1$, and hence, \eqref{eq_m-step_repr} does not
hold. In fact, the proposed methodology applied to the optimal
travel time $T_n^{Opt}$ very soon results in analytically
infeasible calculations. Litvak and van Zwet~\cite{litvak04}
analysed the optimal route. They employed the results on the
$m$-step strategy to derive a recursive expression for the
distribution of the minimal travel time.

We would like to also note that the process of comparing the
lengths of the spacings and deriving corresponding linear
combinations of normalised exponentials can be easily translated
into a computer program. Then, for moderate values of $n$ the
exact distribution of the optimal travel time can be obtained
numerically. The result will be a complicated mixture of linear
combinations of uniform spacings. For large values of $n$ such
exact calculations will require too much computer capacity.
However, in this case, the knowledge of the exact distribution is
not very important since one can apply approximations based on
asymptotic results discussed in the next section.

\subsection{Asymptotic results}\label{ss:asymptotic results}

When the order is large, we can model this situation by letting
$n\to\infty$. Then the expressions in
\eqref{eq_ni_repr} and \eqref{eq_m-step_repr} for the travel time
allow us to obtain asymptotic results that are of independent
mathematical interest. Obviously, if $n\to\infty$ then the travel
time under any strategy goes to one with probability~1. However,
with linear scaling, we obtain non-trivial distributions that we
present below for the nearest-item heuristic and for the optimal
travel time.
\begin{theorem}
\label{theorem2} Let $X_1$, $X_2$, $\ldots$,$X'_1$, $X'_2$, $\ldots$,  be independent
exponentials with mean 1. Then
\begin{align}
\label{eq_ni_limit_distr}
&(n+1)\left(1-T^{NI}_{n}\right)\stackrel{d}{\longrightarrow}
\sum_{j=1}^{\infty}\frac{1}{2^{j-1}}{X_j}\qquad \mbox{\rm (Litvak and Adan~\cite{litvak02})},\\
\label{eq_opt_limit_distr}
&(n+1)\left(1-T^{Opt}_{n}\right)\stackrel{d}{\longrightarrow}
\max\left\{\sum_{j=1}^{\infty}\frac{1}{2^{j}-1}{X_j},
\sum_{j=1}^{\infty}\frac{1}{2^{j}-1}{X'_j}\right\}\quad \mbox{\rm
(Litvak and van Zwet~\cite{litvak04})}
\end{align}
as $n\to\infty$.
\end{theorem}
Result \eqref{eq_opt_limit_distr} is also generalised to the case
when items positions are independent and have some positive
density $f$~\cite{litvak05}.

The expression in the right-hand side of \eqref{eq_ni_limit_distr}
is a well-known functional of the Poisson process, which has been
extensively studied in the literature. We will briefly discuss
this topic in Section~\ref{ss:exponential_functionals}.

\section{Multiple carousels: modelling challenges}\label{s:maria}
The problems examined so far relate to one-carousel models. In
industry though, one rarely meets a facility where only one
carousel is used. Multiple-carousel systems tend to have a higher
level of throughput; however, they increase the investment cost
due to the extra driving and control mechanisms
\cite{hwang91,hwang99}. A natural question is how much the
throughput of a standard carousel can be improved by the
corresponding multiple-carousel system that has the same number of
shelves as the standard carousel. Thus, the question we would like
to examine in this section is the following: given a setup, i.e.\
a specific storage scheme of the items stored on the carousel
\emph{and} a specific travelling strategy, such as those
described in the previous section, how much can we increase the
utilisation of the picker (by assigning to him more carousels to
handle) without increasing the response time of an order above some
chosen level? In other words, how do we reach a \emph{quality and
efficiency} regime in a real situation?

To illustrate things better, consider the following simple
example. A facility assigns $n$ carousels to a single picker. Each
carousel is assigned to an order of a single customer, and each
order consists of exactly one item. Moreover, each carousel
rotates independently until the desired item reaches the picker, who is standing at a fixed point, the origin.
Once this position is reached, the carousel stops until the item
is picked. Only then will the next order be given to the carousel,
which will start rotating the new order to the origin. The picker
serves the carousels in a fixed order, visiting each carousel only
once in every cycle. Clearly, as $n$ goes to infinity, the
utilisation of the picker in steady state tends to one, since
almost surely he will never have to wait. The carousels will have
brought each of their respective items to the origin by the time
the picker is ready to serve them. On the other hand, the time
until the picker returns to the first carousel tends to infinity;
i.e.\ each individual customer suffers long waiting times.

Multiple carousel problems differ intrinsically from
single-carousel problems in a number of ways. Such systems tend to
be more complicated. The system cannot be viewed as a number of
independently operating carousels (cf.~\cite{mcginnis86} and
Section~\ref{ss:Design issues}), since there may be some
interaction between two separate carousels by means of the picker
that is assigned to them. Namely, if the number of pickers is less
than the number of carousels, then the picking strategy that is
chosen for an isolated carousel may affect significantly the
waiting time of another carousel. Thus, one
cannot guarantee that minimising the travel time of a single
carousel maximises the total throughput of the system; the outcome may be quite the
contrary because of the system's interdependency. Another point is
that in multiple-carousel problems, the i.i.d.\ assumption of the time needed to pick each of two consecutive orders with random item storage is in
principle invalid. Characteristics such as the time needed to
reach the optimal point or the travel time for each carousel
depend on one another through the picker's movements. For all
these reasons, multiple-carousel systems merit a special
reference.

Ideally, the problems of minimising the travel time of all
carousels and maximising the picker's utilisation without
surpassing certain levels of each order's response time should be
studied together. However, the interdependence that appears in
multiple carousel problems usually leads to complicated
mathematical structures that can hardly be analysed exactly. One
will have to resort to simplifications.

One technique that can help overcome some of these difficulties is
the setting proposed in Vlasiou {\em et al.}~\cite{vlasiou04}. The
system we consider below consists of two carousels operated by a
single picker. Given a setting, i.e.\ a storage scheme and a
travel strategy, one first needs to obtain an estimate of the
travel time needed in order to collect all items under this
setting. For example, if the items are stored in random positions
on the carousel, then the distribution of the travel time under
the nearest-item heuristic is given by \eqref{eq_ni_dist}. In most
settings though, this distribution cannot be computed
analytically, in which cases the empirical distribution or
simulation may provide a partial answer. Subsequently, one may need
to approximate this distribution by a phase-type distribution; see
e.g.~\cite{osogami}. Then, the following modelling assumption is
made. We aggregate all items in one. That is, we consider an order
that consists of exactly one item. It is assumed that the travel
time of the carousel until that single item is reached is
uniformly distributed (i.e.\ it is assumed that the item is
located randomly on the carousel), while the distribution of the
pick time for that item is taken to be equal to the phase-type
distribution computed previously. Under these assumptions, one can
compute the utilisation of the picker by applying the results
developed in Vlasiou {\em et al.}~\cite{vlasiou04}. This procedure
can be repeated until the desired quality and efficiency regime is
reached.

To describe things concretely, we consider a system consisting of two identical carousels and one picker. At each carousel there is an infinite supply of pick orders that need to be processed. The picker alternates between the two carousels, picking one order at a time. There are two ways one can view this. Either, as mentioned above, one aggregates all items in an order in one super-item (i.e.\ we consider an order that consists of exactly one item) or under the term ``picking time'' we understand the total time needed for the actual picking and travelling from the moment the picker is about to pick the first item in an order until the time the last item is picked. For ease of presentation, we will opt for the first solution, considering orders consisting of exactly one item.

As in Section~\ref{s:nelly}, we model a carousel as a circle of length $1$ and we assume that it rotates in one direction at a constant speed. The picking process may be visualised as follows. When the picker is about to pick an item at one of the carousels, he may have to wait until the item is rotated in front of him. In the meantime, the other carousel rotates towards the position of the next item. After completion of the first pick the carousel is instantaneously replenished and the picker turns to the other carousel, where he may have to wait again, and so on. Let the random variables $A_n$, $B_n$ and $W_n$, $n \ge 1$, denote the pick time, rotation time and waiting time for the $n$-th item. Clearly, the waiting times $W_n$ satisfy the recursion
\begin{equation}\label{recursion}
W_{n+1}=\max\{0,B_{n+1}-A_n-W_n\}, \qquad n = 0, 1, \ldots
\end{equation}
where $A_0 = W_0 \stackrel{\mathsf{def}}{=}0$. We assume that both $\{A_n\}$ and $\{B_n\}$, $n \ge 1$, are sequences of independent identically distributed  random variables, also independent of each other. The pick times $A_n$ follow a phase-type distribution and the rotation times $B_n$ are uniformly distributed on $[0,1)$ (which means that the items are randomly located on the carousels). Then $\{W_n\}$ is a Markov chain, with state space $[0,1)$. Moreover, it can be shown that $\{W_n\}$ is an aperiodic, recurrent Harris chain, which possesses a unique equilibrium distribution.
In equilibrium, equation \eqref{recursion} becomes
\begin{equation} \label{lindley}
W\stackrel{d}{=}\max\{0,B-A-W\}.
\end{equation}
Once the distribution of $W$ is computed from \eqref{lindley}, we can compute $\e[W]$ and thus also the throughput of the system $\tau$ from
\begin{equation}
\label{tau}
\tau = \frac{1}{\e[W]+\e[A]}.
\end{equation}
Equation~\eqref{lindley} with a plus sign instead of minus sign in front of $W$ at the right-hand side, is precisely Lindley's equation for the stationary waiting time in a PH/U/1 queue. The equation for the standard PH/U/1 queue has no simple solution, but in Vlasiou \emph{et al}.~\cite{vlasiou04} we show that the waiting time of the picker in our problem can be solved for explicitly.

For example, assume that the service times follow an Erlang distribution with scale parameter $\lambda$ and $n$ stages; that is,
$$
\Dfa(x)=1-\mathrm{e}^{-\lambda x}\sum_{i=0}^{n-1}\frac{(\lambda x)^i}{i!}, \qquad x \ge 0
$$
and define $\pi_0=\p[W=0]$. Then, for the Laplace transform $\ltw(s)$ of $W$, i.e.\
\[
\ltw(s)=\int_0^1 e^{-sx}\dfw(x)dx,
\]
where $\dfw(x)$ is the density of $W$, the following theorem holds (recall that the domain of integration is bounded by the length of the carousel).
\begin{theorem}[Vlasiou \emph{et al}.~\cite{vlasiou04}]\label{3th:Erlang/U}
For all $s$, the transform $\ltw(s)$ satisfies
\begin{equation}\label{3eq:erl_polywn}
\ltw(s) R(s) = -\mathrm{e}^{-s}s(\lambda+s)^{n}T(-s)-\lambda^{n}T(s),
\end{equation}
where
\begin{align*}
R(s)&=s^2(\lambda^2-s^2)^{n}+\lambda^{2n},\\
T(s)&=\pi_0\biggl(\lambda^{n}+\mathrm{e}^{-(\lambda+s)}\sum_{i=0}^{n-1}\sum_{j=0}^i\frac{s\lambda^i(\lambda+s)^{n-i-1+j}}{j!}\biggr) -\mathrm{e}^{-s}(\lambda+s)^{n}+\\
&\quad+\mathrm{e}^{-(\lambda+s)}\sum_{i=0}^{n-1}\sum_{j=0}^i\sum_{\ell=0}^j\binom{j}{\ell}\frac{s\lambda^i(\lambda+s)^{n-i-1+j}}{j!} \,\ltw[{(\ell)}](-\lambda).
\end{align*}
\end{theorem}
In \eqref{3eq:erl_polywn} we still need to determine the $n+1$ unknowns $\pi_0$ and $\ltw[{(\ell)}](-\lambda)$ for $\ell=0,\ldots,n-1$. Note that for any zero of the polynomial $R$, the left-hand side of \eqref{3eq:erl_polywn} vanishes (since $\ltw$ is analytic everywhere). This implies that the right-hand side should also vanish. Hence, the zeros of $R$ provide the equations necessary to determine the unknowns. In \cite{vlasiou04} it is explained how to determine these unknown parameters (which incidentally form the unique solution to a linear system of equations) and how to invert the transform. Qualitatively, the result is as follows.
\begin{theorem}[Vlasiou \emph{et al}.~\cite{vlasiou04}] \label{3th:densityold}
The density of $W$ on $[0,1]$ is given by
\begin{equation}\label{3eq:density}
\dfw(x)=\suml_{i=1}^{2n+2}c_i \mathrm{e}^{r_i x}, \qquad 0 \le x \le 1,
\end{equation}
and
\begin{equation}\label{3eq:pi0old}
\pi_0 = \p[W=0] = 1-\sum_{i=1}^{2n+2} \frac{c_i}{r_i}(\mathrm{e}^{r_i}-1),
\end{equation}
where $r_i$ is a zero of the polynomial $R$ appearing in Theorem~\ref{3eq:erl_polywn}, and where the coefficients $c_i$ are known explicitly.
\end{theorem}
As a by-product, we have that
\begin{corollary}\label{3cor:tau}
The throughput $\tau$ satisfies
$$
\tau^{-1} = \e[A] + \e[W] =
\frac{n}{\lambda} +
\suml_{i=1}^{2n+2} \frac{c_i}{r^2_i} [1 + (r_i -1)\mathrm{e}^{r_i}] .
$$
\end{corollary}
\begin{remark}
  The same qualitative result holds in case the pick times follow a mixed-Erlang distribution. In this case, the waiting time density is again a mixture of exponentials, where all parameters can be computed explicitly; cf.~\cite{vlasiou04}.
\end{remark}

In a series of papers, Vlasiou \emph{et
al}.~\cite{boxma07,vlasiou07a,vlasiou05,vlasiou05b,vlasiou03,vlasiou04,vlasiou08,vlasiou07}
have relaxed several of the assumptions made above for the
two-carousel setting. For example, the travel time needed to pick all items in an order can have any general distribution (e.g.\ depending on the pick strategy that is followed). In such cases, one
can compute the distribution of the waiting time of the picker by
approximating the distribution of the travel time by an
appropriate phase-type distribution. Phase-type distributions may
be used to approximate any given distribution on $[0,1)$ for the
travel times arbitrarily close; see for example
Asmussen~\cite{asmussen-APQ}. As an illustrative example, we give
below the steady-state distribution of the waiting time of the
picker in case the pick times follow some general distribution
with Laplace-Stieltjes transform (LST) $\lta$, and the travel times follow an
Erlang distribution with parameter $\mu$ and $n$ stages. Here, $\ltw$ denotes the (unknown) LST of the waiting
time of the picker. In this case we have the following:
\begin{theorem}[Vlasiou and Adan~\cite{vlasiou05}]\label{4th:density}
The waiting-time distribution has a mass $\pi_0$ at the origin, which is given by
$$
\pi_0=\p[B<W+A]=1-\sum_{i=0}^{n-1}\frac{(-\mu)^i}{i!}\,\phi^{(i)}(\mu)
$$
and has a density $\dfw$ on $[0, \infty)$ that is given by
\begin{equation}\label{4eq:A density}
   \dfw(x)=\mu^n\mathrm{e}^{-\mu x}\sum_{i=0}^{n-1}\frac{(-1)^i}{i!}\,\phi^{(i)}(\mu)\frac{x^{n-1-i}}{(n-1-i)!}.
\end{equation}
In the above expression, we have that
$$
\phi^{(i)}(\mu)=\sum_{k=0}^i \binom{i}{k} \ltw[(k)](\mu)\,\lta[(i-k)](\mu)
$$
and that the parameters $\ltw[(i)](\mu)$ for $i=0,\ldots,n-1$ are the unique solution to the system of equations
\begin{equation}\label{4eq:balance eq}
\begin{aligned}
\ltw(\mu)&=1-\sum_{i=0}^{n-1}(-\mu)^i\Bigl(1-\frac{1}{2^{n-i}}\Bigr) \sum_{k=0}^i \frac{\ltw[(k)](\mu)\,\lta[(i-k)](\mu)}{k!\,(i-k)!}\\
&\hspace{-2cm}\mbox{and for $\ell=1,\ldots,n-1$}\\
\ltw[(\ell)](\mu)&=\sum_{i=0}^{n-1}\mu^{i-\ell}\frac{(-1)^{i+\ell}}{2^{n-i+\ell}}\frac{(n-i+\ell-1)!}{(n-i-1)!}\sum_{k=0}^i \frac{\ltw[(k)](\mu)\,\lta[(i-k)](\mu)}{k!\,(i-k)!}.
\end{aligned}
\end{equation}
\end{theorem}

As a final curiosity, we present Figure~\ref{fig:throughput}. For single-server queuing models it is well-known that the mean waiting time depends (approximately linearly) on the squared coefficients of variation of the interarrival (and service) times; see also Section~\ref{sss:lindley} for connections of this model to the classical single-server queue.  The results in Figure~\ref{fig:throughput}, however, show that for this two-carousel model, the throughput $\tau$, and thus the mean waiting time, is not very sensitive to the squared coefficient of variation of the pick time; it indeed decreases as $c^2_A$ increases, but very slowly. This phenomenon may be explained by the fact that the waiting time of the server is bounded by one, that is, the time needed for a full rotation of the carousel.

We refrain from giving all results derived for the waiting time distribution in this setting, as they can be found in the papers mentioned so far. One point needs to be stressed though. This technique makes usage of several simplifications (e.g.\ aggregating orders in one item) and approximations (e.g.\ modelling various distributions as a phase-type distribution). Some of them are almost unavoidable. For example, a carousel storing items in separate drawers should be evidently modelled with a discrete travel-time distribution; for the application of these results though, one should approximate this distribution by a (continuous) phase-type distribution. However, the effect that some of these assumptions have to the final result is marginal, or at least fully controlled. As was shown in Vlasiou and Adan~\cite{vlasiou05b}, the error made in computing the distribution of the time the picker has to wait (is not utilised) is bounded. \piccaption{The throughput is almost insensitive to $c^2_A$.\label{fig:throughput}}
\parpic[l]{\includegraphics[width=0.6\textwidth]{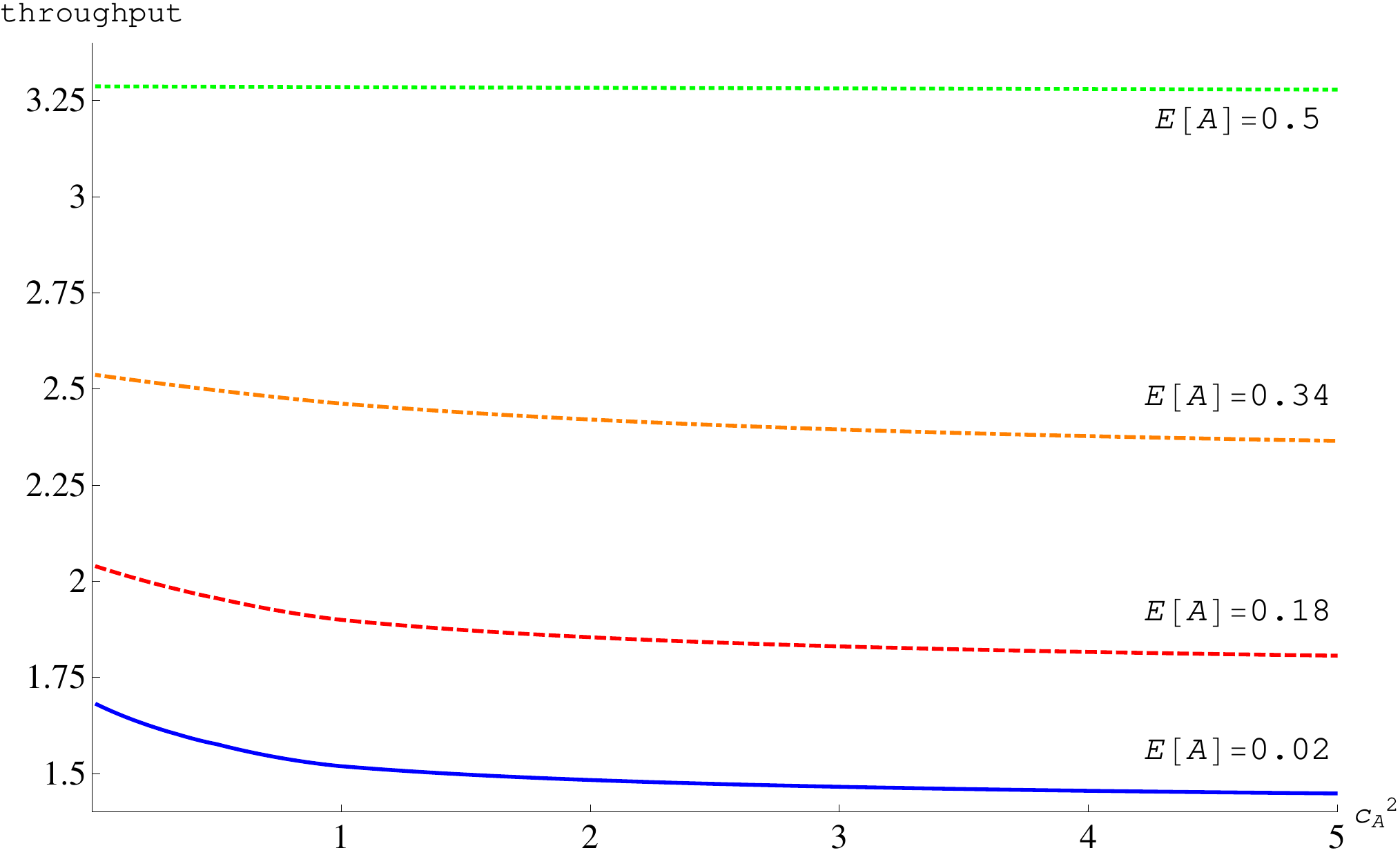}}
%\begin{figure}[hb]
%\begin{center}
%\includegraphics[width=0.8\textwidth]{throughput}
%\end{center}
%\caption{The throughput is almost insensitive to $c^2_A$}
%\label{fig:throughput}
%\end{figure}
Error bounds have been studied widely. The main
question is to define an upper bound of the distance between the
distribution in question and its approximation, that depends on
the distance between the governing distributions.

For our model, recall that $A$, $B$, and $W$ denote respectively the pick
time needed for an item, the travel time of the carousel until
this item is reached, and the waiting time of the picker until the
carousel stops for the pick. Moreover, $\Dfb$ represents the
distribution of $B$ (and similarly also for $W$) and $\eDfb$ is
its approximation (such as the phase-type approximation mentioned above). Using
this approximation, $\eDfb$, one can derive analytically an exact
solution that is obtained for this case for the distribution of
$W$. Denote this solution by $\eDfw$. Then the following error
bound holds.

\begin{theorem}[Vlasiou and Adan~\cite{vlasiou05b}]\label{6th:error bound}
Let $\|\Dfb-\eDfb\|=\varepsilon$. Then $\|\Dfw-\eDfw\|\le {\varepsilon}/(1-\p[B>A])$.
\end{theorem}
In the theorem above, the norm under consideration is the uniform
norm. The main ingredient of the proof relies on the fact that the density for the stationary waiting time of \eqref{recursion} can be described in terms of an integral equation that is a contraction mapping. As a result, approximation errors can be bounded.

An almost identical result can be derived in case one
approximates the pick time, rather than the travel time. Thus, as
this theorem indicates, resorting to approximations yields results
of validity that can be controlled, provided that one has an
estimation of the error that is being made by the original
approximation.

Other results derived for the two-carousel setting include the study of the conditions under which there exists a steady-state distribution~\cite{vlasiou07a}, the study of the tail behaviour of this distribution under general assumptions for the pick and travel times~\cite{vlasiou07a}, the derivation of the steady-state distribution for various cases for the distributions of the pick and travel times~\cite{vlasiou07a,vlasiou05,vlasiou04}, as well as the time-dependent distribution of the waiting times of the picker for a specific setting for the distributions of the pick and travel times~\cite{vlasiou07}. Moreover, certain types of dependencies between the pick and travel times have also been studied, and the steady-state distribution has been derived for these cases as well~\cite{vlasiou08}.

It is worth a mention that such multiple-carousel systems, their
mathematical peculiarities, their results and the way those are
derived are not limited only to carousel, warehousing, or
manufacturing problems. The same equation describing the dynamics
of a two-carousel setting describes also the dynamics of a queuing
model with two nodes that is applied to situations varying from a
university canteen to a surgeon's operating room. For a
description of such systems and detailed analysis see Vlasiou {\em
et al.}~\cite{boxma07,vlasiou07a,vlasiou05,vlasiou07}.

What we have discussed so far on multiple-carousel problems is
summarised as follows. Multiple-carousel problems are
intrinsically different from their single-carousel counterparts.
What is of interest in such problems is to strike a balance
between the utilisation of the picker and the response time of an
order. To date, not much is known about such systems; see
Section~\ref{ss:Problems involving multiple carousels} for an
exhaustive literature review. A few of these results are
simulation studies. However, it is almost inevitable to make use
of some simulation or approximations in these problems. The
results developed in Vlasiou {\em et
al.}~\cite{vlasiou05b,vlasiou04} help predict the performance of
two-carousel systems and ultimately, combined with the results on
e.g nearest-item heuristic or $m$-step strategies discussed in
Section~\ref{s:nelly}, they help design a facility having a
specific quality and efficiency target. However, such results are
still far from accurate. More research is needed on the subject;
specific directions are provided in the next section.

\section{Related research areas}
The mathematics and models involved in the research regarding
carousel systems have surprisingly many connections to broader
areas in queuing theory and applied probability. Other than the
relation to polling systems which will be explained in detail in
Section~\ref{Polling systems}, the subjects we have presented so
far are connected to the classical single-server queue, to
rendezvous networks and layered queues and even to graph theory.
In the following, we highlight few of these connections.

\subsection{Uniform spacings}
\label{ss:spacings}

The uniform spacings defined in \eqref{eq_spacings} constitute a classical mathematical construction which is very well studied.
Uniform spacings have been analysed extensively in two classical
review papers by Pyke~\cite{pyke65,pyke72}. In particular, \cite{pyke65} discusses the connections between uniform
spacings and exponential random variables that are a main concept in the
methodology presented in Section~\ref{s:nelly}. The Markovian property (which is also called the memoryless property) of the exponential distribution is
systematically exploited in Operations Research and in particular
in queuing theory~\cite{asmussen-APQ}.

Uniform spacings play an important role in mathematical
statistics. Mainly, they are used for goodness-of-fit tests which
examine how well a sample of data agrees with a given
distribution $F_0$ as its population. The idea of using uniform
spacings is based on the integral transformation $x\to F_0(x)$
which reduces the problem to the testing of uniformity of the
transformed sample. There is a vast literature on the
distributions, limiting behaviour, approximations and bounds for
various goodness-of-fit test statistics and empirical processes
based on uniform spacings. These investigations are of great
mathematical and practical interest. Considerable progress in the
area has been achieved in the eighties, but there are still many
open problems motivating new studies.

In his detailed review, Pyke~\cite{pyke65} distinguishes two main
types of goodness-of-fit statistics based on a function of uniform
spacings, namely a sum of the form
\[G_n=\sum_{i=1}^ng_n(D_i),\]
or a function of the ordered spacings and their ranks. The
analysis of the first kind of tests goes back to
Le~Cam~\cite{lecam58} and gives rise to an extensive literature, see
e.g.\ \cite{ghosh01,pyke72,wells93} and references therein. Recent
progress on multivariate spacings has been reported in
\cite{li08}. The second type of tests requires the knowledge of
the properties of ordered spacings. This subject has been
extensively studied; we refer the interested reader to the work by
Deheuvels~\cite{deheuvels82} and
Devroye~\cite{devroye81,devroye82}. An original discrete version
of the problem is analysed by Henze~\cite{henze95} who derives the
distribution of the maximal and minimal spacings in lottery
tickets.

Apart from the tests mentioned above, there are also tests based
on $m$-spacings which are the gaps between the order statistics
$U_{i:n}$ and $U_{i+m:n}$. For the analysis of such test
statistics and their asymptotic properties as the number of
observations goes to infinity, see, e.g., Del
Pino~\cite{delpino79}, Hall~\cite{hall86}, and references therein.
The tests based on ordered $m$-spacings have been also analysed,
see, e.g.,~\cite{beirlant85,deheuvels84}. More references on this
subject and results on the approximations for $m$-spacings can be
found in \cite{glaz94}. For further analysis and applications of various
empirical processes based on spacings see Pyke~\cite{pyke72},
Beirlant et al.~\cite{beirlant85,beirlant91}, Einmahl and Van
Zuijlen~\cite{einmahl88,einmahl92} and references therein.

\subsection{Exponential functionals of Poisson processes}
\label{ss:exponential_functionals}

Let $X_1,X_2,\ldots$ be i.i.d.\ exponential random variables with
mean~1. For any $q\in(0,1)$, define
\[J^{(q)}=(q^{-1}-1)\sum_{j=1}^{\infty}(q^{-j}-1)^{-1}X_j,\]
\[I^{(q)}=\sum_{j=1}^{\infty}q^{j-1}X_j.\]
Note that the right-hand side of \eqref{eq_ni_limit_distr} is
exactly $I^{(q)}$ with $q=1/2$. Likewise, the right-hand side of
\eqref{eq_opt_limit_distr} is the minimum of two independent
random variables distributed as $J^{(1/2)}$. We see that the sums
of independent exponentials with exponentially decreasing
coefficients play an important role in the limiting results for
the travel time in carousel systems as the number of items goes to
infinity. Specifically, these random variables appear if we
consider the difference between the travel time and one complete
carousel rotation, and then scale this quantity linearly with the
number of items.

Now let $N(t)$ be a standard Poisson process. Then we can write
$I^{(q)}$ as an exponential functional associated with $N(t)$:
\[I^{(q)}=\int_0^{\infty}q^{N(t)}\,dt.\]
The functional $I^{(q)}$ has been intensively studied in the
literature. Its density was obtained independently in
\cite{bertoin04,dumas02}, and in \cite{litvak01a} for $q=1/2$.
Carmona~\emph{et al.}~\cite{carmona97} derived a density of
$\int_0^{\infty}h(N(t))\,dt$ for a large class of functions $h:
\mathbb{N}\longrightarrow\mathbb{R}_+$, in particular, for
$h(n)=q^n$. Bertoin and Yor~\cite{bertoin04} found the fractional
moments of $I^{(q)}$. If $i^{(q)}(t)$ is a density of $I^{(q)}$,
then $i^{(q)}(t)$ and all its derivatives equal 0 at point $t=0$.
This implies that all moments of $1/I^{(q)}$ are finite. However,
for $q=1/e$, it was proved in \cite{bertoin02} that $1/I^{(1/e)}$
is not determined by its moments. Guillemin \emph{et
al.}~\cite{guillemin04} found the distribution and the fractional
moments of the exponential functional
\begin{equation}
\label{xi} I(\xi)=\int_0^{\infty}e^{-\xi(t)}\,dt,
\end{equation}
where $(\xi(t), t\ge 0)$ is a compound Poisson process.

The distribution function of $I^{(q)}$ and $J^{(q)}$ has an
interesting asymptotic behaviour in the neighbourhood of zero.
Bertoin and Yor~\cite{bertoin02} obtained the following
logarithmic asymptotics:
\[\log i(t)\sim -\frac{1}{2}(\log(1/t))^2\quad {\rm as}\quad t\to +0, \]
where $i(t)$ is a density of
\[I=\int_0^{\infty}e^{-N(t)}\,dt=\sum_{j=1}^{\infty}e^{-j}X_j.\]
The exact asymptotic behaviour  has been derived by Litvak and van
Zwet~\cite{litvak04}. Compared to the logarithmic asymptotics,
their formula contains several additional terms and reveals an
unexpected oscillating behaviour involving theta-functions. The
explanation of why the oscillations appear seems to lie in the
sort of a `binary tree structure' of the functional $I$, whose
coefficients are negative powers of $e$. Later on,
Robert~\cite{robert2005} and  Mohamed and
Robert~\cite{mohamed2005} found that such oscillating asymptotic
behaviour is a typical feature of algorithms with a tree
structure. This phenomenon is compelling and
deserves further studies.

Exponential functionals of Poisson process and, more generally, of
{L}\'evy processes, appear in a number important applications. For
instance, they are relevant to the analysis of randomised
algorithms~\cite{flajolet04} and in mathematical
finance~\cite{bertoin05}. In \cite{dumas02} and \cite{guillemin04}
the exponential functionals of Poisson processes, and,
respectively, of compound Poisson processes, play a key role in the
analysis of the limiting behaviour of a Transmission Control
Protocol connection for the Internet. We refer to the survey
\cite{bertoin05} for further applications, results and references.
The study of exponential functionals of {L}\'evy processes are a current subject of research, see e.g. \cite{lindner09},
\cite{patie09}.

\subsection{Lindley's recursion}\label{sss:lindley}

One of the most intriguing mathematical observations that arise when studying the two-carousel model presented in Section~\ref{s:maria} is that Recursion~\eqref{recursion} differs from the original Lindley's recursion~\cite{lindley52}, which is $W_{n+1}=\max\{0, B_{n} - A_n + W_n\}$, only in the change of a plus sign into a minus sign. At the right-hand side of these two recursions, the sign in front of $W_n$ is reversed. Lindley's recursion describes the waiting time $W_{n+1}$ of a customer in a single-server queue in terms of the waiting time of the previous customer, his or her service time $B_{n}$, and the interarrival time $A_n$ between them. It is one of the fundamental and most well-studied equations in queuing theory. For a detailed study of Lindley's equation we refer to Asmussen~\cite{asmussen-APQ}, Cohen~\cite{cohen-SSQ}, and the references therein.

In the applied probability literature there has been a considerable amount of interest in generalisations of Lindley's recursion, namely the class of Markov chains, which are described by the recursion $W_{n+1}=g(W_n,X_n)$. The model in Section~\ref{s:maria} is a special case of this general recursion and it is obtained by taking $g(w,x)=\max\{0, x-w\}$. Many structural properties of the recursion $W_{n+1}=g(W_n,X_n)$ have been derived. For example Asmussen and Sigman~\cite{asmussen96a} develop a duality theory, relating the steady-state distribution to a ruin probability associated with a risk process. For more references in this domain, see for example Borovkov~\cite{borovkov-ESSP} and Kalashnikov~\cite{kalashnikov02}. An important assumption which is often made in these studies is that the function $g(w,x)$ is non-decreasing in its main argument $w$. For example, in \cite{asmussen96a} this assumption is crucial for their duality theory to hold. Clearly, in the special case of $g(w,x)=\max\{0, x-w\}$ which is discussed in Section~\ref{s:maria}, this assumption does not hold. This fact produces some surprising results when analysing the equation.

The implications of this `minor' difference in sign are rather far reaching. For example, in Section~\ref{s:maria} we have presented two results in Theorems~\ref{3th:Erlang/U} and \ref{3th:densityold}, where we have seen that the waiting time of the picker can be solved for explicitly. For Lindley's recursion, i.e.\ with a plus sign instead of minus sign for $W$ in stationarity, this case correspond to the stationary waiting time in a classical single-server PH/U/1 queue. However, this equation has no simple solution for Lindley's recursion, while we have derived a closed-form expression for the carousel recursion. Furthermore, numerical results (see also Figure~\ref{fig:throughput}) show that for this carousel model the mean waiting time is not very sensitive to the coefficient of variation of the pick time, which is in complete contrast to Lindley's recursion. For these reasons, we believe that it is interesting to investigate in detail the impact on the analysis of such a `slight' modification to the original equation. In this section, we highlight some of the differences of these two models.

\subsubsection{Stability}
For the single-server queue, i.e.\ Lindley's recursion, it is well-known \cite[Ch.\ III.6]{asmussen-APQ} that the random variables representing waiting times of customers converge in distribution (and in total variation) when the mean of the associated random walk is less than zero, or equivalently when the traffic intensity $\rho$ is less than 1; i.e., when $\e[B]<\e[A]$, where we recall that $B$ is the generic service-time random variable, and $A$ is the generic interarrival-time random variable.

For the two-carousel model, though, which is given by Recursion~\eqref{recursion}, the situation is slightly different. In case $\p[B<A]>0$, the stochastic process $ \{W_n\} $ is an aperiodic, (possibly delayed) regenerative process with the time points where $W_n=0$ being the regeneration points. Moreover the process has a finite mean cycle length. To see this, let $X_n=B_n-A_{n-1}$, define the stopping time $\tau=\inf\{n \geqslant 1: X_{n+1}\leqslant 0\}$, and observe that a generic cycle length is stochastically bounded by $\tau$ and that
\[
\p[\tau>n]\leqslant\p[X_k > 0 \mbox{ for all } k=2,\dotsc,n+1]=\p[X_2 > 0]^{n}.
\]
Moreover, we have that $\p[X_2 > 0]<1$ because of the condition $\p[B<A]>0 \Leftrightarrow \p[X<0]>0$. Therefore, from the standard theory on regenerative processes it follows that the limiting distribution exists and the process converges to it in total variation. Through coupling, stability can be shown also for the case where $\p[X<0]=0$; see \cite{vlasiou07a} for details. We see thus that while for Lindley's recursion the stability condition is given by $\e[X]<0$, for Recursion~\eqref{recursion} \emph{stability always holds}; moreover, excluding the deterministic case, we have convergence in total variation.

\subsubsection{Tail behaviour}
For Lindley's recursion, there has been a substantial amount of investigations on the behaviour of $\p[W>x]$ as $x\rightarrow\infty$, the state of the art can be found in \cite{korshunov97}. Results of this type for Recursion~\eqref{recursion} have been derived in \cite{vlasiou07a}. If the right tail of $e^{X}$ is regularly varying of index $-\gamma$ (see \cite{bingham-RV} for background), then
\[
P(W>x) \sim \e[e^{-\gamma W}] \p[X>x].
\]
If the right tail of $e^{X}$ is of rapid variation (see again \cite{bingham-RV}), then
\[
P(W>x) \sim \p[W=0] \p[X>x].
\]
In both equations, we use the notational convention $f(x)\sim g(x)$ to denote $f(x)/g(x)\rightarrow 1$ as $x\rightarrow \infty$.
Note that the class of distributions covering these results include all phase-type distributions, as well as the Weibull, Gamma, Lognormal and Pareto distributions. Moreover, these results indicate that large values of $W$ are caused by a single large value of $X$. This is contrasting with the qualitative picture for Lindley's recursion, where a large value of $W$ is most likely caused only by a single big jump only in the case where $X$ is heavy-tailed. If $X$ is light-tailed (for example phase type), then a large value of $W$ is the cause of a more intricate event involving a change of measure; see \cite{asmussen-APQ} for background.

A natural question is whether it is possible to unify the results for Lindley's recursion and Recursion~\eqref{recursion}. This is possible by considering a recursion that has a minus before $W_n$ (cf.\ Recursion~\eqref{recursion} too) only with probability $1-p$, $p\in [0,1]$, and has a plus before $W_n$ (i.e.\ equal to Lindley's recursion) with probability $p$. For this recursion, the tail behaviour has been studied in \cite{vlasiou09} under assumptions similar to the ones made in \cite{korshunov97}. To summarise the qualitative picture emerging from that paper, the tail behaviour for the unified recursion with $p\in [0,1]$ converges continuously to the results for Recursion~\eqref{recursion} (i.e.\ if $p=0$) for the heavy-tailed case, while it has a discontinuity for $p=1$; for the so-called Cram\'er case the result is reversed: the unified recursion is continuous for $p=1$ and discontinuous for $p=0$, while for the intermediate case (where $X$ is light tailed but does not satisfy the Cram\'er condition) the results for the unified recursion are continuous at both end-points.

\subsubsection{Time-dependent behaviour}
It is well known that for Lindley's recursion, the time-dependent waiting-time distribution is determined by the solution of a Wiener-Hopf problem, see for example \cite{asmussen-APQ} and \cite{cohen-SSQ}. Recursion~\eqref{recursion} though, regularly gives rise to {\em generalised} Wiener-Hopf equation. For example, in \cite{vlasiou07a} we have derived a generalised Wiener-Hopf equation for the density of the stationary waiting time, while \cite{vlasiou07} contains an integral equation for the generating function of the distribution of $W_n$ that is equivalent to a generalised Wiener-Hopf equation, which cannot be solved in general. In Noble~\cite{noble-MBWHT} it is shown that such equations can sometimes be solved, but a general solution, as is possible for the classical Wiener-Hopf problem (arising in Lindley's recursion), seems to be absent.

This makes it appear that \eqref{recursion} may have a more complicated time-dependent behaviour than Lindley's recursion. However, a point we make in \cite{vlasiou07} is that this is not necessarily the case. Thus, Equation \eqref{recursion} is a rare example of a stochastic model which allows for an \emph{explicit} time-dependent analysis. The reason is that, if $B_1$ has a phase-type distribution, we can completely describe \eqref{recursion} in terms of the evolution of a finite-state Markov chain.

We shall refrain from giving all results on the time-dependent behaviour of \eqref{recursion} or their differences from the classical Lindley recursion for the single-server queue, as these results have been well documented elsewhere \cite{vlasiou07a}. Here, we simply list the major findings.

Other than deriving the time-dependent waiting time distribution for \eqref{recursion} under the assumption that the random variables $B_i$ are phase-type distributed, one can derive explicit expressions for the correlation between two waiting times. It results that the covariance function $c(k)$ between two waiting times with lag $k$ converges to zero geometrically fast in $k$. This is consistent with the fact that the distribution of $W_n$ converges geometrically fast to that of $W$, cf.\ Vlasiou~\cite{vlasiou05a}. One of the properties of $c(k)$ is that it is non-negative if $k$ is even and non-positive if $k$ is odd. If in addition, the random variable $X=B-A$ has a strictly positive density on an arbitrary interval, then the inequalities are strict. In contrast, the literature on the covariance function of the waiting times for the single-server queue seems to be sporadic. For the G/G/1 queue, Daley~\cite{daley68} and Blomqvist~\cite{blomqvist68,blomqvist69} give some general properties. In particular, in \cite{daley68} it is shown that the serial correlation coefficients of a stationary sequence of waiting times are non-negative and decrease monotonically to zero.

As we have mentioned before, $\{W_n\}$, as given by \eqref{recursion}, is a regenerative process; regeneration occurs at times when $W_n=0$. Other transient results relate to  the length of a generic regeneration cycle $C$. For Recursion~\eqref{recursion}, we do not need to resort to the usage of generating functions, as is necessary when analysing the corresponding quantity in Lindley's recursion. Note that the interpretation of $C$ for the carousel model is completely different from the corresponding quantity for Lindley's recursion. There, $C$ represents the number of customers that arrived during a busy period. In the carousel setting, $C$ represents the number of pauses a picker has until he needs to pick two consecutive orders without any pause. In this sense $C$ can be seen as a ``non-busy period''.

\subsection{The machine repair problem}
When deriving Equation~\eqref{recursion}, one of the main assumptions we have made, which led to this particular form for the equation is that the picker is not allowed to pick two consecutive orders at the same carousel and must alternate between the two carousels (thus picking all odd-numbered orders from one carousel and all even-numbered orders from the other). This condition is crucial. If we remove this condition, then under certain distributional assumptions, the problem turns out to be the classical machine repair problem, and certain analogies between these two models arise.

In the machine repair problem, there is a number of machines working in parallel (two in our situation, corresponding to the two carousels) and one repairman (corresponding to the picker), who serves the machines when they fail. The machines are working independently and as soon as a machine fails, it joins a queue formed in front of the repairman where it is served in order of arrival. A machine that is repaired is assumed to be as good as new. The machine repair problem, also known as the computer terminal model (see for example Bertsekas and Gallager ~\cite{bertsekas-DN}) or as the time sharing system (see, e.g., Asmussen~\cite[p.\ 79]{asmussen-APQ} or Kleinrock ~\cite[Section 4.11]{kleinrock-QS2}) is a well studied problem in the literature. It is one of the key models to describe problems with a finite input population. A fairly extensive analysis of the machine repair problem can be found in Tak\'{a}cs ~\cite[Chapter 5]{takacs-ITQ}. In \cite{vlasiou05} we compare the two models and discuss their performance.

The issue that is usually investigated in the machine repair problem is the waiting time of a machine until it becomes again operational. In the situation described in Section \ref{s:maria} though, we are concerned with the waiting time of the repairman. It is quite surprising that although the machine repair problem under general assumptions is thoroughly treated in the literature, this question remains unanswered. In the machine repair problem the operating time of the machine is usually more valuable than the utilisation of the repairman, which might explain why the classical literature has been mainly focused on performance measures related to the machines.

In \cite{vlasiou05} the waiting time of the repairman is derived under the assumption that `rotation' times follow a phase-type distribution while `pick' times are generally distributed. Moreover, it is shown that the random variables for the waiting time for the picker/repairman in the two models \emph{are not stochastically ordered}. However, on average, the alternating strategy connected to the two-carousel model leads to longer waiting times for the picker, which readily implies that the \emph{throughput of the machine repair model is bigger}. Furthermore it is shown that the probability that the picker does \emph{not} have to wait is \emph{larger} in the two-carousel alternating system than in the machine repair (i.e.\ non-alternating) model one. This result is perhaps counterintuitive, since the inequality for the mean waiting times of the picker in the two situations is reversed. Regarding the relationship between the $i$-th waiting time of the picker in the two-carousel alternating model (denote this by $W^{\mathrm{A}}_i$), and that of the repairman in the machine repair problem (let this be given by $W^{\mathrm{NA}}_i$), an immediate corollary of the results stated above is as follows.
\begin{corollary}
For all $i$, $\sum_j^i W^{\mathrm{A}}_j \geqslant_{st} \sum_j^i W^{\mathrm{NA}}_j$.
\end{corollary}
So, although the stationary random variables $W^{\mathrm{A}}$ and $W^{\mathrm{NA}}$ are not stochastically ordered, the partial sums of the sequences ${W^{\mathrm{A}}_i}$ and ${W^{\mathrm{NA}}_i}$ are.
Moreover, a conjecture stated in \cite{vlasiou05} suggests that a direct application of the Karlin-Novikoff cut-criterion (cf.\ Szekli~\cite{szekli-SODAP}) leads to an \emph{increasing convex ordering}, namely:
\begin{conjecture}For all increasing convex functions $\phi$, for which the mean exists, we have that
$$\e [\phi(W^{\mathrm{NA}})]\leqslant\e [\phi(W^{\mathrm{A}})].$$
\end{conjecture}

\subsection{Rendezvous networks and layered queues}
The essence of layered queueing (a special case of which is rendezvous networks) is a form of simultaneous resource possession \cite{omari07}.

In its most simple form in computer science applications, in a rendezvous network, a task may serve requests in two rounds (phases) of service. In computer applications, tasks or applications may act both as customers that needs service from other tasks and as servers to other tasks too. As a naive example, think of an application that adds up numbers. It acts both as a server, accepting requests from other applications that need numbers added, and as a customer, requiring service from the central processing unit. One can imagine that tasks are ordered in several levels or layers. Tasks have directed arcs to other tasks at lower layers to represent service requests. A request from an task (client) to a lower-layered task (server) may return a reply to the requester (a synchronous request, or rendezvous). While in the first phase (i.e.\ in the rendezvous) the client is blocked and the server merely continues the thread of control of its client. However, in the second phase the client has an independent thread of control of its own. For example, Task A makes a request to Task B which then makes a request to Task C. While Task C is servicing the request from Task B, Tasks A and B are both blocked \cite{franks09}. Among the advantages of the rendezvous is efficiency, since it provides communication without the effort of buffer management and the message copying associated with asynchronous communication. However, some potential for concurrency is lost, and there may be performance-impairing bottlenecks when a key task spends long periods send-blocked \cite{neilson95}. Special approximations are needed to solve queueing models which contain a two-phase server, because the second phase effectively creates a new customer in the queueing network, violating the conditions of product form queueing \cite{franks09}.

Distributed or parallel software with synchronous communication via rendezvous is found in client-server systems and in proposed Open Distributed Systems, in implementation environments such as Ada, V, Remote Procedure Call systems, in Transputer systems, and in specification techniques such as CSP, CCS and LOTOS. The delays induced by rendezvous can cause serious performance problems, which are not easy to estimate using conventional models which focus on hardware contention, or on a restricted view of the parallelism which ignores implementation constraints. Stochastic Rendezvous Networks are queueing networks of a new type which have been proposed as a modelling framework for these systems. They incorporate the two key phenomena of \emph{included service} and the \emph{second phase of service} mentioned above. The main work on rendezvous networks focuses on Mean Value Analysis and gives approximate performance estimates. This method has been applied to moderately large industrial software systems \cite{woodside95}.

A Layered Queuing Network (LQN) model is a canonical form for extended queueing networks that represent layered service. In a layered queue a server, while executing a service, may request a lower layer service and wait for it to complete. Thus, in LQNs there exist entities that have a dual role; they act as servers to other entities of a lower layer and as customers to higher layered entities. The service time of the upper server includes the queueing delay and service time of the lower server, and this may extend through multiple layers. LQN was developed for modelling software servers, with for example blocking remote procedure calls to lower layer software servers, however it applies to any extended queueing network in which resource usages are nested, lower layer usages within higher layer usages \cite{omari07}.

The two-carousel model we have presented in Section~\ref{s:maria} is a layered queue, and in particular a rendezvous network. To see this, organise the system as follows. The items that are stored on the carousel and have to be picked comprise the lowest layer. Carousels are in the middle layer, while the picker is put in the highest layer. One may view the rotation time of a carousel as a first phase of service for the item that will be picked. The carousel (middle layer) acts in this case as a server. However, the second phase of service (the actual picking) does not necessarily happen immediately (rendezvous). The item might have to wait for the picker to return from the previous carousel -- cf.\ Recursion~\eqref{recursion}. At this stage, the carousels act as customers waiting to be served by the higher layer, the picker. We see thus that each carousel acts both as a server (rotating items to the picking location) and as a customer (waiting until the picker completes his task before the carousel can resume its role as a server, bringing the next item to the picking location).

Layered systems are quite unknown outside the computer-science community. E.g., in \cite{perel08} it is mentioned that ``this paper presents a model, never studied before in the queueing literature, of a system of two connected queues where customers of one queue act as the servers of the other queue'' -- a comment that may very well be valid outside the computer-science literature. The analysis of Recursion~\eqref{recursion}, as it developed in \cite{boxma07,vlasiou07a,vlasiou05,vlasiou05b,vlasiou04,vlasiou09a,vlasiou07} as well as \cite{perel08} are the only papers we are aware of that deal with LQNs using analytic and probabilistic tools, and admittedly all the aforementioned work on the two-carousel model had not made the connection between this model and layered queues.

\subsection{Maximum weight independent sets in sparse random graphs}
However unusual it might be in queuing theory to encounter a non-increasing Lindley-type recursion, Recursion~\eqref{recursion} appears in problems involving the computation of the distribution of the maximum weight of an independent set in a sparse random graph.

Consider an $n$-node sparse random $r$-regular graph (i.e.\ a graph selected uniformly at random from the set of all graphs on $n$ nodes in which every node has degree $r$). An independent set is a set of nodes of the graph where no two nodes in the set are connected by an edge. Suppose that the nodes of the graph are equipped with some nonnegative weights $w_i$ which are generated independently according to some common distribution $F_w$. One may be interested for example in the limits of maximum weight independent sets and matchings in sparse random graphs for some types of i.i.d.\ weight distributions. Then Recursion~\eqref{recursion} corresponds exactly to the one related to the weight distribution in an 1-regular graph; see \cite{gamarnik06}. Moreover, if one considers $r$-regular graphs, then the corresponding recursion giving the weight distribution in this case is similar to the one corresponding to the waiting time of a picker serving $r$ carousels; see \eqref{eq:multi-carousel}. The crucial difference in this case is in \eqref{eq:multi-carousel} the random variables $W_{n+1}$ and $W_n$ appearing at the right-hand side of the recursion are not independent, while the corresponding variables in the recursion related to $r$-regular graphs are independent; see \cite[Eq.\ (3)]{gamarnik06}. It would be interesting to investigate the connections between the research areas of warehouse logistics and graph theory.

\section{Further research}\label{s:research}
\subsection{Considering different item storage schemes}
As mentioned in Section~\ref{s:nelly}, as of yet the case of independent
uniformly distributed items locations is the only known scenario
where the travel time can be evaluated analytically by applying  a systematic
mathematical approach. It is important to develop methods to obtain statistical characteristics of the travel time
under more realistic assumptions on the items locations. As we
discuss below in Section~\ref{ss:Picking a single order}, there
are not many results in this direction in the literature. The
non-uniform distributions of pick positions and especially the
correlations between the items in an order lead to challenging
mathematical problems. We believe that no feasible analytical
solutions can be obtained in most of the realistic models. Thus,
the problem calls for well justified heuristics and efficient
numerical methods.

\subsection{Further topics in two-carousel problems}
The model we have considered in Section~\ref{s:maria} applies to a
two-carousel system that is operated by a single picker.
Two-carousel systems have received some attention in the
literature (cf.\ Section~\ref{ss:Problems involving multiple
carousels}) but many questions remain open. A line of research is
directed towards studying the performance of two-carousel systems
under various storage-assignment policies (randomised or not), for
various pick/travel time strategies and heuristics (sequential
picking, nearest-item heuristic, $m$-step strategies, etc.), for
single- or dual-command cycles, and for open- and closed-loop
strategies. Here a single command cycle assumes a single
operation, such as only storage or only retrieval. In a
dual-command cycle, a storage and retrieval are combined to
efficiently use the time of the operator. Furthermore, an open-loop strategy
implies that the carousel remains stationary at the point where
the last item was retrieved (awaiting the next order to be fed),
while under the closed-loop strategy the carousel returns to a
predefined point after the retrieval of an order is completed.

As
explained in Section~\ref{s:maria}, two-carousel systems differ in
nature and in analysis from the corresponding one-carousel
problems even when studied under the same assumptions on the
various storage, pick, cycle, and starting-point strategies that
are followed. Since two-carousel systems perform in broad terms
better than single-carousel systems \cite{hwang99}, studying the
expected increase of the throughput of the system can help answer
questions of financial nature, such as whether the benefits from
the increased throughput justify the increased cost of building
and operating a two-carousel system.

\subsection{Extensions to multiple carousels}
The model discussed in Section~\ref{s:maria} can be extended to
the case of multiple carousels  as follows. For instance,
consider the situation where a single picker operates
three carousels. Apart from the number of carousels, all
other characteristics of the model remain the same as in
Section~\ref{s:maria}. That is, we consider again an infinite
queue of orders that need to be picked, we have again
a rotation stage and a picking stage for
each item. Moreover, as before, the picker serves all carousels cyclically. For three carousels, this leads to the recursion
\begin{equation}\label{eq:multi-carousel}
W_{n+2}=\max\{0,B_{n+2}-W_{n+1}-A_{n+1}-W_{n}-A_n\},
\end{equation}
where now the variables appearing at the right-hand side are not
independent of one another, as was the case for all variables
appearing at the right-hand side of Recursion~\eqref{recursion}. We may assume for convenience that the
sequences $\{A_n\}$ and $\{B_n\}$ are independent among them and
between them. Furthermore,
we note that the waiting times $W_{n}$ and $W_{n+1}$ are not
independent. The state of the system can be modelled e.g.\ as a
two-dimensional Markov chain, where apart from the waiting time of
the picker for the $n$-th item that will be picked we also need to incorporate the
remaining rotation time of the next carousel to be served.
Evidently, if the rotation times are assumed to be
exponentially distributed, the system (for three or more carousels)
can be analysed explicitly by similar techniques as the ones
applied in Chapter~4 of~\cite{vlasiou}, although it is doubtful how realistic such an assumption is.

Naturally, if one considers a system with multiple carousels or
stations, one can think about optimisation questions. Namely, as
the number of carousels increases, the waiting time of the picker
is expected to decrease. After serving a long series of carousels
cyclically, when you return to the beginning of the cycle, with
high probability the item to be picked will have reached the
origin. This implies that an item will have to wait for the picker
at the origin more frequently than in the two-carousel system,
which means that the throughput of a single carousel decreases.
Intuitively, as the number of carousels increases to infinity, the
utilisation of the picker increases to one, while the throughput
of each individual carousel decreases to zero. Given a setting,
one might wonder how many carousels a single picker can operate so
that we maximise both the throughput of the carousels and the
utilisation of the picker simultaneously.

\subsection{Incorporating picking strategies to multiple carousel problems}
The ultimate goal of the analysis of  carousel systems is to
provide a mathematical model that adequately describes the reality
and, at the same time, can be efficiently evaluated either
analytically or numerically. At the moment, the literature on a
single carousel has advanced enough to characterise the travel
time with great precision, at least for independent uniform items locations. However, as mentioned above, single
carousel systems are rarely used in modern warehouses. Clearly,
multiple carousel models are more relevant from a practical point of
view. The drawback is that such models tend to become extremely
complex. Until now the studies of multiple carousel systems were
either solely based on simulations or employed analytical models
that involved simplifying assumption on the order picking time.
For instance, in Section~\ref{s:maria} we assumed that each order
is collected within a random time that has the same distribution
for each order. This is definitely a simplifying assumption, because,
for instance, the orders may differ in size, and as we saw
in Section~\ref{s:nelly}, the distribution of the travel time
depends on the number of items to be collected. Further literature
on multiple carousels discussed in detail in
Section~\ref{ss:Problems involving multiple carousels} also
involves significant simplifications of the real-life situation.

In this respect, a major challenge for future studies is to
develop a unified approach for rigorous studies of real-life
automated  storage and retrieval systems. Such an approach is expected to
 involve the methods proposed so far for single and
multiple carousels. In Sections~\ref{s:nelly}~and~\ref{s:maria} we
presented well-developed methodologies for analytical studies of
order picking in one and two carousel units. Thus, an important
topic for further research is to combine these two problems in one
integrated study of multiple carousel systems. One may hope to
obtain interesting analytical results in this direction because of
the analytical nature of both methodologies. However, the problems
of combining these two settings are challenging. In Section
~\ref{s:nelly} we have seen that the travel time distribution can
be of a complicated form, while the results in
Section~\ref{s:maria} often rely on assumption such as exponential
or phase-type pick times (recall that the travel time needed to
pick all items corresponds to the pick time for orders aggregated
in one item). Also, as mentioned above, the travel time depends on
the size of the order, while the technique of aggregating orders
in one item has made use of the assumed independence between pick
times and rotation times (while one might expect that in orders
with multiple items, long travel times might be correlated to
orders with multiple items and thus to shorter rotation times to
the first item in the order). Eventually, one will have to resort
to the development of reasonable algorithms rather than the
derivation of exact distributions. In this respect, we emphasise
again that algorithmic studies of realistic carousel models
constitute an important part of further research.

\subsection{Considering the order arrival process}
It is also interesting to study if single or multiple carousel systems can
be analysed in case there is an arrival process according to which
the orders arrive. If orders arrive according to a Poisson process in front of the carousel,
what can be said for the waiting time of the
picker? This question can also be combined with a
non-alternating system, where the picker serves the first carousel
that has completed the rotation to the next item on that carousel that needs to be picked, or with Bernoulli-type
requests, where the picker has to serve with a certain probability
the ``first'' carousel and with the complementary probability
the ``other'' carousel (potentially waiting for an item if none
is present at the designated carousel). For each case, one should
also consider the stability of the system in case the arrival rate
of the orders is less than the throughput of the system with an
infinite queue of orders.

\subsection{Polling systems}\label{Polling systems}
A polling system is comprised of a number of customer queues that are served in an order by a single server. In the literature on polling systems, the polling system with two
queues where at each queue the server serves exactly one customer
before switching to the other queue is often referred to as the
\textit{1-limited alternating-service} model. The model described in Section~\ref{s:maria} is closely related to such polling systems. The two main differences are the existence of an extra stage, the rotation time of the carousel, and the absence of an arrival process for the orders. In polling systems one deals only with one stage, which in the terminology of Section~\ref{s:maria} is represented by the picking stage. Extending the model
of Section~\ref{s:maria} by introducing an arrival process of the
orders as suggested above, is equivalent to studying an
1-limited alternating-service model with switch-over times between
the stations (which can be seen as being equivalent to the
rotation time towards the single item).

The polling model with two queues,
Poisson arrivals, and no switch-over times has first been studied
by Eisenberg~\cite{eisenberg79}, where the main question studied
is the
queue-length distribution, as is often the case in the literature on polling systems. Eisenberg~\cite{eisenberg79} gives the
generating function for the stationary joint distribution of the
two queue sizes. Cohen and Boxma~\cite{cohen81} study the single
server queue with two Poissonian arrival streams and no
switch-over times. The server handles alternatingly a customer of
each queue if the queues are not empty and it is assumed that
customers of the same arrival stream have the same service time
distribution. It is shown that the determination of the joint
queue-length distribution at the departure epoch can be formulated
as a Riemann-Hilbert boundary problem that can be completely
solved for general service time distributions. Introducing
switch-over times increases the complexity of the problem. In
Boxma~\cite{boxma85} the analysis is extended to include
switch-over times of the server between queues, under the
restriction that both queues have identical characteristics. This
work is further extended in Boxma and Groenendijk~\cite{boxma88},
where the authors no longer request that both queues have
identical characteristics. It is assumed that service times and
switch-over times are generally distributed.

The literature on polling systems with alternating service is not
limited to the references above but is rather extensive; see
\cite{groenendijk,ibe90,ozawa90} for some references. It seems
though, that the question regarding the waiting time of the picker
for the 1-limited polling system with two
carousels has not been considered outside the scope
of~\cite{vlasiou}. Thus, introducing an arrival process for the
orders in the model of Section~\ref{s:maria} complements the
existing literature on polling systems and forms a challenging
problem. The interesting feature then is that the switch-over time
between two queues depends on the current picking time. Again, the
results from Section~\ref{s:nelly} can be incorporated into the
model for adequate description of order picking times.

An extension considered in polling systems is the
\textit{$k$-limited} service policy, where the server switches
queues after having served at most $k$ customers in one queue. For
an extensive list of references on $k$-limited polling systems see
Van Vuuren and Winands~\cite{vuuren06}. The main focus of the
existing literature is again on the queue-length distribution of
all stations. As the authors note in \cite{vuuren06}, ``to this
very day, not only hardly any exact results for polling systems
with the $k$-limited service policy have been obtained, but also
their derivations give little hope for extensions to more
realistic systems''. It is worth considering the $k$-limited
service discipline under the exact setting we have established in
Section~\ref{s:maria}, where now the focus is on the distribution
of the waiting time of the server.

\section{Literature overview}\label{s:literature}
In the following, we classify the literature on carousels
according to the main theme handled. This taxonomy allows for a
better overview of the variety of the subjects examined. A crucial
distinction is made between systems that involve a single carousel
and systems with multiple carousels. The first four categories
presented relate to single-carousel systems, while systems with
multiple carousels are examined later on.

\subsection{Storage}\label{ss:Storage}
The performance of a carousel system depends greatly upon the way
it is loaded and the demand frequency of the items placed on it.
An effective storage scheme may reduce significantly the travel
time of the carousel. Several strategies have been followed in
practice to store items on a carousel. The simplest strategy is to
place the items randomly on the carousel. \textit{Randomised
policies} have been examined extensively \cite{hwang91,litvak},
and various performance characteristics have been derived under
the assumption that the items are uniformly distributed on the
carousel.

One way to improve the throughput of a carousel system is to adopt
a storage policy other than the randomised assignment policy. Ha
and Hwang~\cite{ha94} have studied what they call the
``two-class-based storage'', which is a storage scheme that
divides the items in two classes based on their demand frequency.
The items with a higher turnover are randomly assigned to one
continuous region of the carousel, while the less frequently asked
items occupy the complementary region. The authors show by
simulation that the two-class-based storage can reduce
significantly the expected cycle time, both in the case where a
cycle is a single pick or storage of an item (single-command
cycle), and in the case where a cycle consists of the paired
operations of storing and retrieving (dual-command cycle). The
same authors in \cite{hwang94} examine the effects of the
two-class-based storage policy on the throughput of the system,
and present a case where there is a 16.29\% improvement of this
policy over the randomised policy.

Another storage scheme is suggested by Stern~\cite{stern86}.
Assignments are made using a \textit{maximal adjacency principle},
that is, two items are placed closely if their probability of
appearing in the same order is high. The author evaluates this
storage assignment analytically by using a Markov chain model he
develops.

The \textit{organ pipe arrangement} for a carousel system is
introduced in Lim {\em et al.}~\cite{lim85} and is proven to be
optimal in Beng\"{u}~\cite{bengu95} and in Vickson and
Fujimoto~\cite{vickson96} under a wide variety of settings. The
organ pipe arrangement has been widely used in storage units, such
as magnetic tapes \cite{bitner79} and warehouses
\cite{malmborg90}. This arrangement is based on the classical
mathematical work of Hardy, Littlewood and Polya~\cite{hardy-I}.
Their concept is used in \cite{bitner79} to minimise the expected distance travelled
by an access head as it travels from one record to
another. Various optimality properties of this arrangement have
been proven; see for example Keane {\em et al.}~\cite{keane84} and
references therein.

\piccaption{Illustration of the organ pipe arrangement, where the upper numbers indicate the frequency ranking of an item.\label{test}}
\parpic[l]{\includegraphics[width=0.4\textwidth]{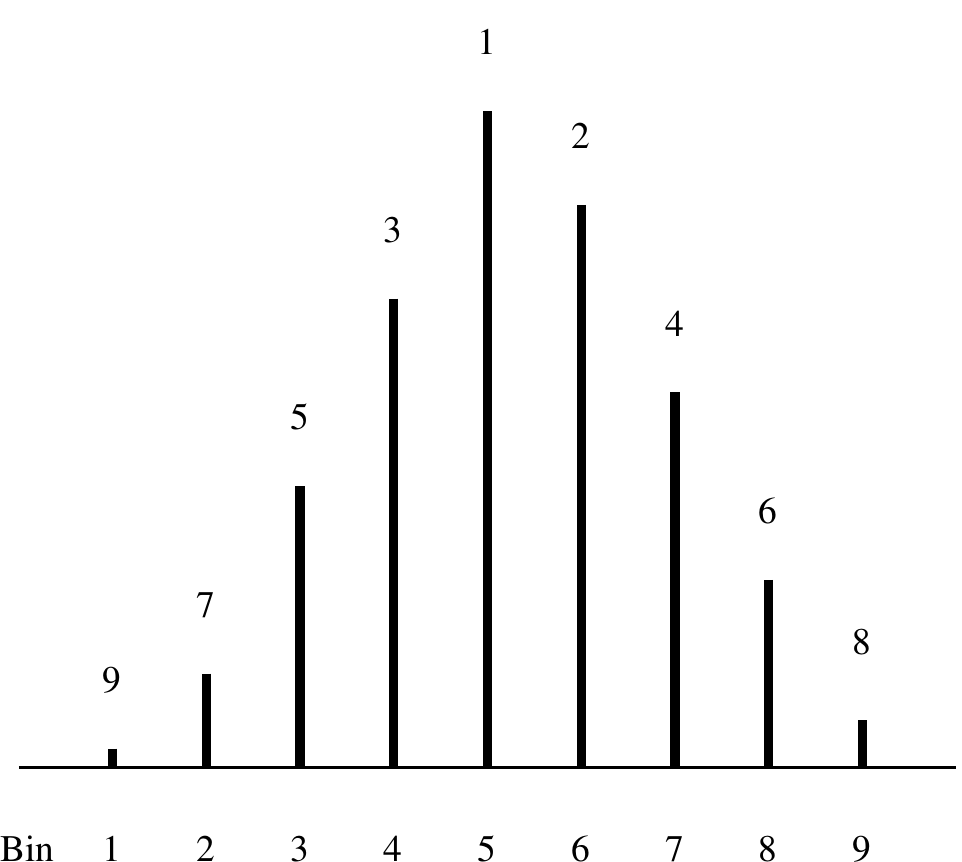}}
In carousel systems, the organ pipe arrangement places the item with the
highest demand in an arbitrary bin, the items with the second and
third highest demands in the bin next to the first one but from
opposite sides, and sequentially all other items next to the
previous ones, where the odd-numbered items according to their
frequency are grouped together and placed next to one another in a
decreasing order from the one side of the most frequent item (and
similarly the even-numbered items are grouped together and placed
to the other side). Figure~\ref{test} illustrates the organ pipe arrangement.
The numbers at the top indicate the ranking of an item in a decreasing order of frequency.
%\begin{figure}[htbp]
%\leavevmode
%\begin{center}
%\includegraphics[width=0.6\textwidth]{organpipe}
%\end{center}
%\caption{Illustration of the organ pipe arrangement, where the upper numbers indicate the frequency ranking of an item.}
%\label{arrangement}
%\end{figure}

Park and Rhee~\cite{park05} study the system throughput and the
job sojourn times under the organ pipe arrangement, where
independent one-item orders arrive according to a Poisson process.
They explicitly quantified the gain of the organ pipe arrangement
compared to random assignment and showed that this gain grows with
the  `skewness' in the demand distribution.

Abdel-Malek and Tang~\cite{abdel94} study the travel times in
carousels with $N$ bins and the organ pipe arrangement under the
assumption that each order consists of one item and a sequence of
orders forms a Markov chain: if the current order requires bin $p$
then the next order requires bin $q$ with probability $P_{pq}$.
The objective is to find the optimal assignment, which minimises
the average travel time. Their extensive numerical experiments
show that although the organ-pipe arrangement is not optimal in
this setting, it performs very close to optimality in a wide range
of system parameters. The optimal solution in~\cite{abdel94} is
determined by solving a {\it quadratic assignment problem}. The
quadratic assignment problem is a well-known optimisation problem
on choosing an optimal permutation of $n$ coordinates of a vector
${\bf x}=(x_1,\ldots,x_n)$ in order to minimise ${\bf x}C{\bf
x}^T$, where $C$ is a cost matrix. Such problems have a long
history started with the work of Koopmans and
Beckmann~\cite{koopmans57}. Litvak~\cite{litvak05} shows by
experimental studies and by providing asymptotic results for large
orders that in general the optimal storage depends on the order
size. Moreover, the organ-pipe storage is disadvantageous when an
order is large.

Another question related to storage is about the number of items
of each type that should be stored on the carousel in order to
maximise the number of orders that can be retrieved without having
to reload. This question is examined in Jacobs {\em et
al.}~\cite{jacobs00}, where the authors propose a heuristic that
yields a reasonable solution, the error of which can be bounded.
This method has been improved by Yeh~\cite{yeh02}, where a more
accurate solution is obtained, and further on by Kim~\cite{kim05},
where it is observed that the heuristic described in \cite{yeh02}
does not always lead to the optimal solution. The author
constructs an algorithm that yields the optimal solution. This
algorithm is further improved in Li and Wan~\cite{li05}. This line
of research has been continued in the recent paper by
Hassini~\cite{hassini08}. In the formulation used in Jacobs {\em
et al.}~\cite{jacobs00}, the author determines the optimal
allocation. Along with exact optimal solutions for deterministic
and stochastic demand, \cite{hassini08} also provides heuristics
that perform close to optimal.

\subsection{Picking a single order}\label{ss:Picking a single order}
One of the most important performance characteristics of a
carousel system is the total time to pick an order. The total time
to retrieve all items of an order may be expressed as a sum of the
total time that the carousel is travelling plus the total time
that the carousel is stopped for picking. The latter is
effectively the total pick time, and it is not affected by the
sequence in which we choose to retrieve the objects. However, the
total travelling time greatly depends upon the retrieval sequence.
The analysis of the travel time under various strategies is, in
general, a non-trivial problem. This problem,
however, has been resolved for independent and uniformly distributed item
locations~\cite{litvak}, as we discussed in detail in Section~\ref{s:nelly}.

Various picking strategies have been proposed. Bartholdi and
Platzman~\cite{bartholdi86} assume a discrete model and study the
performance of an algorithm and three heuristics that determine an
efficient, but not necessarily optimal, sequence of retrieving all
items. A heuristic is a simpler, non-optimal procedure that is
based on a specific strategy. The heuristic methods proposed are
the \textit{nearest-item heuristic}, where the next item to be
picked is always the one that is closer to the picker at any given
moment, the \textit{shorter-direction heuristic}, where the
carousel chooses the shortest direction between the route that
simply rotates clockwise and the route that rotates
counterclockwise, and the
\textit{monomaniacal heuristic}, that always chooses to rotate to
the right and pick items sequentially. The \textit{optimal
retrieval} algorithm that is presented enumerates all possible
paths; therefore, it is guaranteed to find the quickest sequence
in which to retrieve a single order.

In \cite{bartholdi86} the authors prove among other things that
the travel time under the nearest-item heuristic is never greater
than one rotation of the carousel. Litvak {\em et
al.}~\cite{litvak01a} provide the upper bound of $1-1/2^n$ full
rotations, where $n$ is the number of items in the order, and show
that the new upper bound is tight. Litvak and Adan~\cite{litvak01}
obtained the distribution and the asymptotic properties of the
travel time under the nearest-item heuristic for uniformly
distributed independent items locations. These results, based on
properties of uniform spacings, have been discussed in detail in
Sections~\ref{ss:narest-item},~\ref{ss:asymptotic results}. In
\cite{litvak01a}, the first two moments of the travel time and the
distribution of the number of turns are computed recursively by
conditioning on the event that there is an empty space of size $x$
on one side of the picker's current position. We presume that such
methods may lead to the travel time distribution in some special
cases with non-uniform items locations.

Another interesting picking strategy that has been already
discussed in Section~\ref{ss:m-step} is the so-called
\textit{$m$-step strategy}, where the carousel chooses the
shortest route among the ones that change direction at most once,
and only do so after collecting at most $m$ items. In case of
independent uniformly distributed items locations the average
travel time under the $m$-step strategy is smaller than the one
under the nearest-item heuristic already for
$m=2$; see \cite{litvak02}. The results by Litvak and
Adan~\cite{litvak02} on the $m$-step strategies have been
presented in Section~\ref{ss:m-step}. In an earlier paper,
Rouwenhorst \emph{et al.}~\cite{rouwenhorst96} apply analytical
methods to study the case when $m\le 2$. This means that the
carousel changes direction after collecting at most two items.
They interpret $m$-step strategies as stochastic upper bounds for
the minimal travel time and present convincing numerical results
on the excellent performance of such strategies.

Wen and Chang~\cite{wen88} model the carousel as a discrete
bidirectional loop and assume that the time to move between the
bins of a shelf is \emph{not} negligible. They propose three
heuristic solution procedures and compare their performance. An
earlier version of this work can be found in Wen~\cite{wen86}.

Ghosh and Wells~\cite{ghosh92} model the carousel as a continuum
of \textit{clusters} and \textit{gaps}, where a cluster is a
segment on the circle that corresponds to a series of locations
that have to be visited for the retrieval of an order, while a gap
is the segment of the circle between two clusters. The authors
develop two algorithms to find optimal retrieval strategies. In particular, to find an optimal path,
they avoid a complete enumeration by noticing that a turn can never be made after covering more than 1/3 length of the carousel.

Stern~\cite{stern86} studies properties of the optimal, i.e.\
minimal, picking sequence both for the \textit{open-loop}
strategy, where the carousel remains stationary at the point where
the last item was retrieved (awaiting the next order to be fed),
and for the \textit{closed-loop} strategy, where the carousel
returns to a predefined point after the retrieval of an order is
completed. He formally shows that under the open-loop strategy the
carousel will change its direction at most once when following the
optimal picking sequence, while under the closed-loop strategy the
carousel will turn at most twice. A recursive expression for the
distribution of the minimal travel time needed to collect one
order of $n$ randomly distributed items in the open-loop scenario
is given explicitly by Litvak and Van Zwet~\cite{litvak04}.

The case when positions of the items in an order are dependent has
not received much studies. One way to model the dependencies is
described by Abdel-Malek and Tang~\cite{abdel94} who assume that
the positions of successive items form a Markov chain. In this
setting, they study the performance of the organ-pipe storage
rule. Stern~\cite{stern86} introduces correlations between items
in an order by considering several order types, where each type
corresponds to a fixed list of items. The work of Wan and
Wolff~\cite{wan04} focuses on minimising the travel time for
``clumpy'' orders, that is, orders concentrated on a relatively
small segment of the carousel, and introduces the nearest-endpoint
heuristic for which they obtain conditions for it to be optimal.
In this setting, one can no longer assume that the items locations
are uniformly distributed. Moreover, there is clearly a strong
dependence between items positions.

The model with non-uniform items locations reflects a relevant
situation when some of the items are required more frequently than
others. Most of the papers that assume distinct frequencies assume
the orders of one item (see e.g.~\cite{bengu95}). An interesting
work on non-uniformly distributed items is given by
Litvak~\cite{litvak05}, where the focus is on the length of the
shortest rotation time needed to collect a single order when the
order size is large and the items locations have a non-uniform
continuous distribution with a positive density $f$ on $[0,1]$.

\subsection{Picking multiple orders}\label{ss:Picking multiple orders}
A popular strategy for reducing the mean travel time per order in
carousel storage and retrieval systems is batching together a
number of orders and then picking them sequentially. A
\textit{batch} is a set of orders that is picked in a single tour.
Two consecutively picked items do not necessarily belong to the
same order. An excellent literature survey by Van den
Berg~\cite{vdberg99} on planning and control of warehousing
systems addresses this issue and the problems that arise if large
batches are formed. Apart from the questions mentioned before,
Stern~\cite{stern86} also considers the performance of a carousel
for a fixed set of order types (for example, big orders with many
items, and small ones).

Bartholdi and Platzman~\cite{bartholdi86} are mainly concerned
with sequencing batches of requests in a bidirectional carousel.
They specify the number of orders to be retrieved (ignoring any
new arrivals) and propose three heuristic methods to solve this
static problem. Orders may be picked in any sequence (and not
necessarily at the order they arrive), but picks within the same
order are performed consecutively. They define the \textit{minimum
spanning interval}, which is the shortest interval containing all
the items of an order and, by assuming that the picker always
begins and finishes retrieving an order at one of the endpoints of
this interval, they construct the shortest matching chain by
ordering the orders accordingly. This procedure may fail to give
an uninterrupted sequence in which to pick the orders; therefore,
they propose the following heuristics. The first one, called the
\textit{hierarchical} heuristic, picks any order that happens to
have a common endpoint with another order, and then travels
clockwise until an unpicked endpoint is encountered, and repeats
the procedure. The \textit{nearest-order} heuristic is practically
an extension of the nearest-item heuristic described earlier in
the paper, as is the case with the \textit{second monomaniacal}
heuristic they propose. Under these heuristics, they obtain upper
bounds for the travel time.

Ghosh and Wells~\cite{ghosh92} assume that the orders have to be
picked under a FIFO sequencing restriction, which means that the
first order to arrive at the warehouse is the first order that
will be picked, and so on. Since the orders are retrieved in a
FIFO fashion, the problem is reduced to finding how to retrieve
each individual order so that the best overall retrieval is
achieved. They develop an algorithm for the optimal retrieval path
of $n$ orders via dynamic programming, and show how to update
dynamically the solution when new orders arrive.

Rouwenhorst {\em et al.}~\cite{rouwenhorst96} model the carousel
as an M/G/1 queuing system, where the orders are the
``customers'' that require service, and the service they get
depends on the pick strategy that is followed. This approach
permits the derivation of various queuing characteristics such as
the mean response time and the waiting time when orders arrive
randomly. The authors mention that the tight upper bounds for the
mean response time can be further exploited to obtain also good
approximations for excess probabilities of the response time.

Van den Berg~\cite{vdberg96} assumes either a fixed or an
arbitrary sequence of orders. When the sequence of the orders is
given, he presents an efficient dynamic programming algorithm that
finds an optimum path that visits all orders in the specified
sequence. Furthermore, when there is no given order sequence, he
simplifies the problem to a \textit{rural postman problem} on a
circle and solves this problem to optimality. The rural postman
problem is the problem of finding the shortest route in an
undirected graph which includes all edges at least one time. Van
den Berg~\cite{vdberg96} concludes that the obtained solution
requires at most 1.5 revolutions more than a lower bound of an
optimal solution to the original problem. Simulation results
suggest that the average rotation time may be reduced up to 25\%
when allowing a free order sequence. Lee and Kuo~\cite{lee08}
formulate the problem of optimal sequencing of items and orders as
a {\it multi-travelling salesman problem}. In the multi-travelling
salesman problem, there are several salesmen in a home city, and
each of the other cities has to be visited only by one salesman.
Using this formulation, Lee and Kuo~\cite{lee08} provide
efficient heuristics for optimal picking of several orders
consisting of multiple items.

\subsection{Design issues}\label{ss:Design issues}
All research papers mentioned so far that deal with travel time models of carousel systems assume average uniform velocity of the carousel. In other words, the main assumption is that the carousel travels with constant speed and the acceleration from the stationary position (when a pick is performed) to its full speed, as well as the deceleration from the maximum speed to zero speed, are negligible factors when computing the travel time of the carousel. Guenov and Raeside~\cite{guenov89} give some empirical evidence that the error induced when neglecting acceleration and deceleration of an order picking vehicle is indeed negligible. Thus the problem of minimising retrieval times can be considered to be equivalent to the problem of minimising the average distance travelled by the carousel per retrieval.

Hwang {\em et al.}~\cite{hwang04}, however, develop strategies for picking that take into consideration the variation in speed of the carousel. For unit-load automated storage and retrieval systems there are several travel-time models that consider the speed profiles of the storage and retrieval robot. In \cite{hwang04} some relevant references are given. Unlike the unit-load automated storage and retrieval systems, almost all the existing travel-time models for carousel systems assume that the effects of the variation in speed are negligible. In \cite{hwang04} the authors try to bridge this gap in the literature. They assume that the items are randomly distributed on the carousel and derive the expected travel time both in the case of a single command cycle and in the case of a dual command cycle. They verify the accuracy of the proposed models by comparing the results to results directly obtained from discrete racks.

Egbelu and Wu~\cite{egbelu98} study the problem of pre-positioning the carousel in anticipation of storage or retrieval requests in order to improve the average response time of the system. Choosing the right starting point of a carousel in anticipation of an order is also referred as the \textit{dwell point selection problem.} This strategy becomes relevant when the items are stored under the organ pipe arrangement. In this situation the dwell point should be chosen to be the location of the most popular item; see, e.g., \cite{bengu95}.

Spee~\cite{spee96} is concerned with developing design criteria for carousels. He states the basic conditions for designing an automatic order picking system with carousels and comments on the optimal storage design. Namely, he is interested in finding the right number of picking robots and the right number and dimensions of a carousel so that the investment is minimised, provided that the size of the orders that need to be retrieved is given.

McGinnis~\cite{mcginnis86} studies some of the design and control issues relevant to \textit{rotary racks}. A rotary rack is an automated storage and retrieval system that strongly resembles carousels. In fact, conceptually, a rotary rack is simply a carousel, where the only difference is that each level or shelf of this carousel can rotate independently of the others. The author concludes that, while rotary racks appear to be a simple generalisation of conventional carousels, the control strategies that have been shown effective for carousels do not appear to be as effective for these systems. Rotary racks can be viewed as a multiple-carousel system (where each level is considered as a sub-carousel) with a single picker.

\subsection{Problems involving multiple carousels}\label{ss:Problems involving multiple carousels}
While almost all work mentioned in this section concerns one-carousel models, real applications have triggered the study of models involving multiple carousels. The study of such models is not as developed yet as the study of models involving a single carousel. The list of references that follows seems to be complete.

Perhaps the first academic study that investigates the performance of a system involving several carousels is that of Emerson and Schmatz~\cite{emerson81}. The authors simulated the operation of the warehouse of Rockwell's Collins Telecommunications Products. The system consists of twenty-two carousels, where each pair of carousels had a single-operator station (so there are in total eleven operator stations). The questions they are concerned with are how big the batch size of orders should be so as to complete the week's work (which is used as a performance measure) and keep all operators busy, what happens when a carousel or a station is down, and how is an overload or an imbalance (for example, unequal operator performance, unequal carousel loading, or large orders) handled. In order to investigate potential solutions to these three imbalance conditions, the authors investigate two operating rules.

The first operating rule studies six different storage schemes
with seven carousel pairs (and thus seven operators). It uses
simulation models to study simple storage schemes such as random
storage, sequential alternating storage, and storage in the
carousel with the largest number of openings. The aim in
\cite{emerson81} is to study the degree of carousel usage. The
authors find that there is no significant difference between the
carousel loads among the storage schemes. However, they do not
treat the problem of optimally assigning items to carousel bins,
and do not present any analytical models to help investigate the
problem. The second operating rule they investigate is a floating
operator. This is an operator who is trained to work at any
station, and who is moving to different stations according to
specific needs (for example, depending on the size of the queue at
a particular station). They conclude that this solution seems
advantageous for the purposes of the warehouse they investigate.

Koenigsberg~\cite{koenigsberg86} presents analytic solutions for
evaluating the performance of a single carousel, and discusses the
ways in which his approach can be extended to a system involving
two unidirectional carousels both served by a single robotic
operator. The carousels are related only through the state of the
robot, which means that each carousel is independent of the other
except for the time it waits for an operation to commence (such as
pick, storage, or repair) because the robot is busy at the other
carousel. The author concludes that under some conditions, it is
often more advantageous to have two carousels of identical length
instead of one carousel of double the length. Furthermore, going
to three carousels of equal length (i.e.\ one third of the length
of the single carousel) will offer little further improvement.

Hwang and Ha~\cite{hwang91} study the throughput performance both
of a single and of a double carousel system. Based on a randomised
storage assignment policy, cycle time models are developed for
single and dual commands. Furthermore, they examine the value of
the information on the succeeding orders in terms of system
efficiency, which may lead to better scheduling of the orders to be
processed.

In a later work, Hwang {\em et al.}~\cite{hwang99} attempt to
measure analytically the effects of double shuttles of the storage
and retrieval machine (i.e.\ the robotic picker) on the throughput
both of the standard and of the double carousel system. Storage
and retrieval machines with double shuttles are machines that have
space for two items. Thus, for example, an item can be retrieved
from the carousel and stored on one shuttle, while the other
shuttle has an item that needs to be stored to the carousel. After
this item is stored, a second item can be retrieved from the
carousel and placed on the empty shuttle. All these operations
occur during a single cycle of the carousel operation. For the
double carousel system, a retrieval sequence rule is proposed
which utilises the characteristics of the two independently
rotating carousels. From the test results, double shuttles are
shown to have a substantial improvement over single shuttles. This
improvement tends to be more prominent in the double carousel
system. Due to cost concerns, the authors note that an economic
evaluation will be needed to justify the extra cost of double
carousel systems and double shuttles before implementing them in
real world situations.

Wen {\em et al.}~\cite{wen89} consider a system comprised of two
carousels and a single retrieval machine. Their main assumption is
that every order must be picked in a single tour, i.e., an order
cannot be divided into two or more sub-tours. Batching orders
together is also not allowed. They analyse the retrieval time and
propose four heuristic algorithms for the scheduling sequence of
retrieving items from the system to satisfy an order. Their method
is an extension of the algorithm presented in \cite{bartholdi86}
and \cite{stern86}.

Meller and Klote study the throughput of a group of several
carousels, a so-called carousel pod~\cite{meller04}. They use
approximations to evaluate the order pick time in one carousel and
then evaluate the throughput of a pod by plugging in the average
response times of each unit and modelling the pod of $c$ carousels
as a queuing system where $1/c$ picker operates one carousel.
Further, they derive an approximation for the system's throughput
using a diffusion approximation by Gelenbe~\cite{gelenbe75} which
was earlier applied by Bozer and White~\cite{bozer96} in the
analysis of end-of-aisle order-picking systems.

Recently, Hassini and Vickson~\cite{hassini03} studied storage
locations for items, aiming to minimise the long-run expected
travel time in a two-carousel setting with a single picker. They
assume that the products are available at all times (so as to be
able to ignore possible delays due to  lack of stock), and that
orders are not batched; that is, the carousel system processes
only single-item orders. This is applicable in situations where
individual product orders are processed in a
first-come-first-served policy, or when the next item to be
retrieved is known only after the present one has been picked. The
authors compare the performance of three heuristic storage schemes
and a genetic algorithm~\cite{goldberg} that for small-sized
problems completely enumerates the solution space. They conclude
that none of the heuristic approaches leads to a solution that
outperforms the algorithmic solution they provide.

The same model is also studied by Park {\em et al.}~\cite{park03}.
As is the case in \cite{hassini03}, in \cite{park03} the basic
assumptions are that  there is an infinite number of items to be
picked and that an order consists of a single item. The authors,
however, are not interested in storage issues. They further assume
that the single operator, the picker, is alternately serving the
two carousels. This may cause the picker to have to wait for an
amount of time until the item at the carousel he is currently
serving is rotated in front of him. They derive the distribution
of the waiting time of the picker under specific assumptions for
the pick times. This allows them to derive expressions for the
system throughput and the picker utilisation.

The model presented in \cite{park03} has been extended further in
Vlasiou {\em et al.}~\cite{vlasiou,vlasiou05b,vlasiou08,vlasiou03,vlasiou04} by
removing all assumptions related to the pick times or rotation
times. In related work, Vlasiou {\em et
al.}~\cite{vlasiou07a,vlasiou05,vlasiou07} have shown that the
two-carousel model studied in \cite{hassini03,park03} is
equivalent to an alternating service queue, if one allows for
rotation times with an infinite support. Some of these results have been presented in Section~\ref{s:maria}.

Finally, we would like to mention that there is a broad literature
on automated storage and retrieval systems (see e.g. the survey by
Le-Duc~\cite{leduc}). An extensive list of references has been
also made available on-line by Roodbergen~\cite{roodbergen_list}.

\section*{Acknowledgements}

The authors would like to thank Ivo Adan for the suggestion to
write this survey paper and for his active involvement in the
research described in Sections~\ref{s:nelly}~and~\ref{s:maria}. In
the Netherlands, the three universities of technology have formed
the 3TU Federation. This article is the result of joint research
in the 3TU Centre of Competence NIRICT (Netherlands Institute for
Research on ICT).

\phantomsection
\addcontentsline{toc}{section}{References}

%\bibliographystyle{apt}
%\bibliographystyle{plain}
%\bibliography{maria}

%\end{document}

\end{document}